\documentclass[11pt,a4paper]{article}
\usepackage[utf8]{inputenc}
\usepackage{graphicx,verbatim,array,multicol,courier}
\usepackage{amsmath,amssymb,amsthm,mathrsfs, mathtools}
\usepackage{color, graphics, graphicx, siunitx}
\usepackage{amsfonts, dsfont, bm}
\usepackage{natbib} 
\usepackage{geometry}
\usepackage{enumitem,float}
\usepackage{lscape}
\usepackage[pdftex,bookmarks,colorlinks]{hyperref}
\usepackage[dvipsnames]{xcolor}
\usepackage{tabularx}
\usepackage[toc,title,page]{appendix}
\usepackage{siunitx}
\usepackage[labelformat=simple]{subfig}
\usepackage[linesnumbered,ruled,vlined]{algorithm2e}
\usepackage{booktabs} 
\usepackage{minitoc}

\newlist{inparaenum}{enumerate}{2}
\setlist[inparaenum,1]{label=(\alph*)}
\setlist[inparaenum,2]{label=(\roman{inparaenumi}\emph{\alph*})}

\makeatletter
\def\adl@drawiv#1#2#3{%
        \hskip.5\tabcolsep
        \xleaders#3{#2.5\@tempdimb #1{1}#2.5\@tempdimb}%
                #2\z@ plus1fil minus1fil\relax
        \hskip.5\tabcolsep}
\newcommand{\cdashlinelr}[1]{%
  \noalign{\vskip\aboverulesep
           \global\let\@dashdrawstore\adl@draw
           \global\let\adl@draw\adl@drawiv}
  \cdashline{#1}
  \noalign{\global\let\adl@draw\@dashdrawstore
           \vskip\belowrulesep}}
\makeatother

\numberwithin{equation}{section}
\theoremstyle{definition}
\newtheorem{defi}{Definition}[section]
\newtheorem{cond}[defi]{Condition}
\theoremstyle{plain}
\newtheorem{theo}[defi]{Theorem}

\newtheorem{lem}[defi]{Lemma}
\newtheorem{cor}[defi]{Corollary}

\theoremstyle{remark}
\newtheorem{rem}[defi]{Remark}

\newcommand{\diff}{\mathrm{d}}

\newcommand{\supp}{\text{supp}}
\newcommand{\ind}{\mathbbm{1}}

\def\ind{ {{\rm 1}\hskip-2.2pt{\rm l}}}
\definecolor{navy}{rgb}{0,0,0.502}
\definecolor{brown}{rgb}{0.59, 0.29, 0.0}
\def\indic{\mathds{1}}

\hypersetup{colorlinks,%
	citecolor=blue,%
	filecolor=green,%
	linkcolor=red,%
	urlcolor=violet,%
}

\newcommand{\Real}{\mathbb{R}}
\newcommand{\Nat}{\mathbb{N}}

\newcommand{\Prob}{\mathbb{P}}
\newcommand{\Expect}{\mathbb{E}}

\newcommand{\bftheta}{{\boldsymbol{\theta}}}
\newcommand{\bfvartheta}{{\boldsymbol{\vartheta}}}

\newcommand{\bfmu}{{\boldsymbol{\mu}}}
\newcommand{\bfbeta}{{\boldsymbol{\beta}}}

\newcommand{\bfSigma}{{\boldsymbol{\Sigma}}}

\newcommand{\bfJ}{{\boldsymbol{J}}}

\newcommand{\bfA}{{\boldsymbol{A}}}

\newcommand{\bfX}{{\boldsymbol{X}}}

\newcommand{\bfS}{{\boldsymbol{S}}}

\newcommand{\bfb}{{\boldsymbol{b}}}
\newcommand{\bfI}{{\boldsymbol{I}}}

\newcommand{\bfm}{{\boldsymbol{m}}}

\newcommand{\lm}{\widetilde{m}}
\newcommand{\lk}{\widetilde{k}}

\DeclareMathOperator*{\argmax}{arg\,max}

\title{Asymptotic theory for the likelihood-based block maxima method in time series}
\author{David L. Carl, Simone A. Padoan and Stefano Rizzelli}

\begin{document}

\maketitle
\begin{abstract}
This paper develops a rigorous asymptotic framework for likelihood-based inference in the Block Maxima (BM) method for stationary time series. While Bayesian inference under the BM approach has been widely studied in the independence setting, no asymptotic theory currently exists for time series. Further results are needed to establish that BM method can be applied with the kind of dependent time series models relevant to applied fields.
To address this gap we first establish a comprehensive likelihood theory for the misspecified Generalized Extreme Value (GEV)  model under serial dependence. Our results include uniform convergence of the empirical log-likelihood process, contraction rates for the Maximum Likelihood Estimator, and a local asymptotically Gaussian expansion. Building on this foundation, we develop the asymptotic theory of Bayesian inference for the GEV parameters, the extremal index, $T$-time-horizon return levels, and extreme quantiles (Value at Risk). Under general conditions on the prior, we prove posterior consistency, $\sqrt{k}$-contraction rates, Bernstein–von Mises theorems, and asymptotic coverage properties for credible intervals. For inference on the extremal index, we propose an adjusted posterior distribution that corrects for poor coverage exhibited by a naive Bayesian approach. Simulations show excellent inferential performances for the proposed methodology.
\end{abstract}

%
\section{Introduction}\label{se:intro}
%

%
\subsection{Frontiers and challenges in extreme value theory}\label{sec:EVT}
%

The aim of Extreme Value Theory (EVT) is to develop statistical models and tools for assessing the risk of future events that are significantly more extreme than those observed historically. This goal is of central importance across a wide range of applied domains, including finance and economics \citep[e.g.,][]{embrechts1997modelling, novak2011extreme}, as well as in climate, weather, and environmental sciences \citep[e.g.,][]{coles2001, beirlant2006, dey2016extreme}. In this work, we focus on the Block Maxima (BM) method, one of the most widely used approaches in the univariate case.

High-frequency data, such as daily financial returns, often exhibit serial dependence that cannot be ignored without risking biased or overly optimistic estimates of extreme events. To address this issue, we consider a strictly stationary time series $(X_t)_{t\geq1}$, in particular with a stationary cumulative
distribution function $F$. Under mild conditions on $F$ and the dependence structure \citep[see, e.g.,][Ch. 3]{leadbetter1983}, the distribution of block maximum $\max(X_1, \ldots, X_m)$ 
can be approximated by a rescaled version of the \textit{generalized extreme value} (GEV) distribution, denoted by $G_{\gamma}^\theta((x - \mu)/\sigma)$, for sufficiently large block size $m$. Here, $\gamma$ is the shape parameter or \textit{extreme value index}, $\mu$ and $\sigma > 0$ are location and scale parameters, respectively, and $\theta \in (0,1]$ is the \textit{extremal index}. The parameter $\gamma$ governs the tail heaviness of the distribution, while $\theta$ reflects the degree of clustering of extreme values, roughly corresponding to the reciprocal of the mean size of clusters of high threshold exceedances \citep[e.g.,][]{leadbetter1983, embrechts1997modelling, beirlant2006}.

This result has several important implications. When the stationary distribution $F$ is unknown and a large sample of size $n$ is observed, the distribution of the block maxima over sufficiently large blocks of size $m$ can be approximated by $G_\gamma^\theta((\cdot - \mu)/\sigma)$. This forms the basis for statistical inference on extremes. Leveraging the GEV approximation, two key risk measures can be inferred:
(1) the $T$-time-horizon return level, defined as the $\tau$-quantile of the block maximum distribution, with $\tau = 1 - 1/T$ for a given time horizon $T$; and
(2) the $\tau$-quantile of the stationary distribution $F$ for extreme levels $\tau \to 1$, often referred to as the \textit{extreme quantile} or \textit{Value at Risk (VaR)}, a commonly used risk metric in finance and economics, \citet[see e.g.][p. 381]{beirlant2006}.

While the BM method is widely employed under the assumption of independence, its theoretical underpinnings are far from straightforward \citep[e.g.,][]{ferreira2015, dombry2019}. Regardless of the inferential method—frequentist or Bayesian—the analysis must account for model misspecification: the GEV distribution only provides an asymptotic approximation, while the block size $m$ must be fixed in practice. Moreover, scale and location parameters depend on $m$, in contrast with classical parametric models where parameters are intrinsic (e.g., \citealp{vaart1998}). Futhermore, the GEV family is statistically irregular, as the sign of $\gamma$ determines its support. Classical regularity conditions (e.g., those of \citealp{cramer1946}) and even weaker Lipschitz-type conditions do not hold \citep{vaart1998, buecher2017, padoan2024}. As a result, asymptotic analysis of likelihood-based inference—including Bayesian approaches—faces substantial challenges. In the time series setting, the presence of extremal clustering further complicates inference, as it reduces the effective amount of information compared to the independence case, increasing the uncertainty in estimating $F$ and its extremal features.

Although the estimation of the extremal index $\theta$ has received considerable attention in the literature \citep[e.g.,][]{smith1994, weissman1998, ferro2003, suveges2007, robert2009, northrop2015efficient}, and the asymptotic theory of its maximum likelihood estimator (MLE) has been developed by \citet{berghaus2018}, no corresponding results are available for Bayesian inference of $\theta$, nor for inference on the GEV parameters $\boldsymbol{\vartheta} = (\gamma, \mu, \sigma)^\top$ in the time series context.

In this paper, we develop the asymptotic theory for both maximum likelihood and Bayesian estimation  
of the GEV parameters $\boldsymbol{\vartheta}$, as well as for key functionals such as $T$-time-horizon return level and extreme VaR \citep[e.g.][Ch. 10]{beirlant2006}. In addition, we provide the first theoretical results for Bayesian inference on the extremal index $\theta$, offering a comprehensive framework for likelihood-based extreme value analysis in dependent time series.

%
%
%
\subsection{Purpose and main results}\label{sec:goals}
%
%
%

Over the past few decades, numerous studies have explored Bayesian inference for the BM method under the assumption of independence \citep[e.g.,][]{coles1996bayesian, coles2001, coles2003anticipating, stephenson2004bayesian, beirlant2006, northrop2016}. However, a rigorous theoretical foundation for Bayesian inference in this setting has only recently been developed by \citet{padoan2024}. To the best of our knowledge, no asymptotic results currently exist for Bayesian inference in the context of BM for time series.

The primary contribution of this article is to establish a rigorous asymptotic theory for Bayesian inference in the BM framework when applied to stationary time series data. To this end, the first part of the paper develops a comprehensive likelihood theory under serial dependence. Specifically, we study the empirical log-likelihood process for the misspecified GEV model, deriving its uniform convergence, convergence of its derivatives (detailed in the supplement), as well as local and global bounds, and a local asymptotically Gaussian expansion. As a relevant byproduct, we also obtain contraction rates for the corresponding MLE. Notably, even in the independence case, whether the MLE—computed over the full parameter space—is consistent and uniquely maximizes the likelihood has remained an open question. We resolve this question for both independent and dependent data, thereby closing a long-standing gap in the literature.

Building on this foundation, the second part of the paper develops the asymptotic theory for Bayesian inference in the time series setting. We provide general and practical conditions on the prior distribution for the GEV parameters under which the corresponding posterior distribution satisfies key asymptotic properties: consistency with a $\sqrt{k}$-contraction rate, a Bernstein–von Mises (BvM) theorem, and asymptotically valid coverage probabilities for credible intervals. These results extend those of \citet{padoan2024} by allowing for more general conditions and by addressing the more complex time series setting. Our framework accommodates informative proper priors for the shape parameter $\gamma$ and data-dependent informative priors for the location and scale parameters. Unlike the results in \citet{padoan2024}, which rely on restrictive assumptions about the density of $F$, our theory is derived under weaker and more standard conditions \citep[see Ch. 2 in][]{haan2006}.
It is worth noting that while a version of the BvM theorem exists for misspecified models \citep{kleijn2012bernstein} 
in
the independence setting, it does not directly apply to our time series setting and the GEV family violates the smoothness and regularity assumptions required by \citet{kleijn2012bernstein}.

We also propose a new Bayesian procedure for inference on the extremal index $\theta$, a quantity central to time series extremes. The standard posterior distributions based on the classical likelihood function used in frequentist inference \citep{berghaus2018} performs poorly in this context, as the corresponding credible intervals fail to achieve the nominal coverage level asymptotically—even in the absence of asymptotic bias. To overcome this problem, we propose an adjusted posterior distribution that yields credible intervals with correct asymptotic coverage, and further show it enjoys usual posterior asymptotic properties.
Finally, given the practical relevance of risk assessment for future extremes, we extend our Bayesian asymptotic theory to include inference on key functionals such as the $T$-time-horizon return level and extreme VaR. We prove consistency, contraction rates, asymptotic normality, and valid asymptotic coverage for the posterior distributions of these quantities. These results provide practitioners with powerful tools for reliably quantifying uncertainty in the magnitude of future extreme events.

%
%
\subsection{Article organization}\label{sec:structure}
%
%
Section \ref{eq:backg} introduces the key concepts, notation, and foundational conditions for GEV-based modeling in the context of time series. Section \ref{sec:LAT} presents the BM approach and develops the asymptotic theory for the MLE of the GEV parameters, as well as for return levels and extreme VaR, by analyzing the properties of the empirical log-likelihood process. Section \ref{sec:Bayes_inf} establishes the asymptotic theory for Bayesian inference in the GEV framework for time series, covering the posterior behavior of the GEV parameters, the extremal index, return levels, and extreme VaR. Section \ref{sec:simulations} provides an extensive simulation study that assesses the finite-sample performance of the proposed methodology. All proofs are collected in the Supplement, which also includes additional technical results, posterior computation details, and further simulation outputs. References to Propositions, Lemmas, etc., that use letters concern results of the Supplement.

%
\section{Statistical framework}\label{eq:backg}
%
We work with a stationary time series $(X_t)_{t\geq 1}$, with stationary cumulative
distribution function $F$. Assume that $F$ is an unknown continuous distribution satisfying the domain of attraction condition, denoted as $F\in\mathcal{D}(G_\gamma)$, with $\gamma\in\Real$ (e.g., \citealp[][Ch. 3]{leadbetter1983}, \citealp[][Ch. 8]{embrechts1997modelling}, \citealp[][Ch. 10]{beirlant2006}).
Under mild conditions on the serial dependence (see Condition D in \citealp[][Chapter 3.2]{leadbetter1983}), namely if the dependence between observations decays sufficiently fast as the time separation between them increases, the domain of attraction property implies in the time series context the following limiting results.

First, like in the independence case, if
\begin{equation}\label{eq:DoA_IID}
\lim_{m\to\infty} F^m(a(m)z + b(m)) = G(z),\quad z\in\Real
\end{equation}
for some sequence sequence $a(m)>0$ and $b(m)\in\Real$, for $m=1,2,\ldots$,  then $G$ must be the GEV distribution 
$$
G_\gamma(z)=\exp\left(-(1+\gamma z)_+^{-1/\gamma}\right),
$$
%
where $(x)_+=\max(0,x)$.  The case $\gamma=0$ reads as $G_{\gamma}(z)=\exp(-\exp(-z))$. 
The extreme value index $\gamma$ describes the weight of the right-tail distribution. $G_\gamma$ is a short-, light- or heavy-tailed distribution whose support is $\mathbb{G}_\gamma=(-\infty, -1/\gamma)$, $\Real$, $(-1/\gamma, \infty)$, depending on whether $\gamma<0$, $\gamma=0$, $\gamma>0$, respectively \citep[][Ch. 1]{haan2006}.  The GEV density is $g_\gamma(z)=G_\gamma(z) \left(1+\gamma z\right)^{-1/\gamma-1}$ and becomes $\exp(-z - \exp(-z))$ in the case $\gamma = 0$. Functions $a$ and $b$ are not unique, however,  the choice $b(t)=Q(\exp(-1/t))$, for $t>1$, is convenient for the derivation of our theoretical results and is then used in the sequel, where, for $\tau \in [0, 1]$, $Q (\tau) = F^{\leftarrow} (\tau) =  \inf( x : F(x) \leq \tau)$ is the generalised quantile function of $F$. Second, let  $M_m=\max(X_1,\ldots,X_m)$ be the sample maximum for $m=1,2,\ldots$ and $\overline{F}=1-F$.
If for any $X>0$ there is a sequence $u_m=u_m(x)$ such that $m\overline{F}(u_m)\to x$ as $m\to\infty$ and 
\begin{equation}\label{eq:DoA_extremal_index}
\lim_{m\to\infty} \Prob(M_m\leq u_m)=\exp(-\theta x),\quad \theta \in(0,1],
\end{equation}
then the time series is said to have an extremal index \citep[see][Theorem 3.7.1, p. 67]{leadbetter1983}. In the sequel we assume this holds.
The extremal index $\theta$ measures the extremal serial dependence by representing the reciprocal of the asymptotic mean size of clusters of consecutive suitably normalised excess variables above an asymptotic threshold.
Because $F \in \mathcal{D}(G_\gamma)$, then one can take $u_m$ of the form $a(m) z+ b(m) $, for $x=(1+\gamma z)^{-1/\gamma}$.
Third, let $(\widetilde{X}_{1},\ldots,\widetilde{X}_{m})$ be iid random variables with marginal distribution $F$, then as $m\to\infty$ we have 
\begin{equation}\label{eq:DoA_TS}
\Prob \left(M_m \leq u_n \right)\approx \Prob \left( \max_{1\leq i\leq \lfloor \theta m \rfloor} \widetilde{X}_i\leq u_n \right)\to G^\theta_\gamma(z),
\end{equation}
see \citet[][Theorem 3.6.6]{leadbetter1983}. This is an insightful result as it suggests that serial dependence impacts the limiting distribution of the normalised sample maximum as if it were computed on a sequence of iid variables of reduced length $\lfloor \theta m \rfloor$, where $\lfloor x \rfloor$ is the largest integer smaller than $x$. Consequently, clustering of large observations implies that there is comparatively less information available to make inference about $F$ in time series than in the independence case and a higher uncertainty is therefore expected in the estimation step.  The above results lead to the following practical implications. 

First, for large $m$, thanks to max-stability of the GEV family \citep[][p. 9]{haan2006}, the distribution of the sample maximum $M_m$ can be approximated by $G_\bfvartheta$, where $\bfvartheta=(\gamma, \mu, \sigma)^\top$, with $\mu$ and $\sigma$ that are representative of the norming constants $b_{\theta}(m)= b(m) + a(m)(\theta^\gamma - 1)/\gamma$ and $a_{\theta}(m)=a(m)\theta^\gamma$, respectively.
The corresponding density is $g_\bfvartheta(x)=g_\gamma((x-\mu)/\sigma)/\sigma$ whose support is $\mathbb{G}_\bfvartheta=(-\infty, \mu-\sigma/\gamma)$, $\Real$, $(\mu-\sigma/\gamma, \infty)$ depending on whether $\gamma<0$, $\gamma=0$, $\gamma>0$. 
The model $G_\bfvartheta$ offers a valuable approximation for the distribution of  $M_m$,  providing a basis for statistical inference when the true distribution $F$ is unknown, as is typically the case in applications. However, for any fixed $m$, this model remains statistically misspecified, regardless of how large $m$ is. Consequently, inferential procedures based on it may suffer from bias. A standard strategy to control the error introduced by such misspecification is to use the second-order condition (\citealp[][Ch. 2 and Appendix B]{haan2006}), which characterizes the rate of convergence of the limit in \eqref{eq:DoA_IID}, thus allowing one to quantify and manage the resulting approximation error. We report it next.
%
\begin{cond}\label{cond:second_order}
The limit \eqref{eq:DoA_IID}, known also as  {\it first-order} condition is equivalent to the following result \citep[e.g.][Theorems 1.1.6 and 1.2.1]{haan2006}. There exists a positive function $a(\cdot)$ such that for all $x>0$
 \begin{equation*}\label{eq:first_order}
\lim_{s\to\infty} \frac{b(sx)-b(s)}{a(s)}=\int_1^{x}u^{\gamma-1}\diff u.
\end{equation*}
Assume that:
\begin{inparaenum}
\item \label{en_sec_order_1} There is a positive function $a(\cdot)$ and a positive or
negative function $A$ satisfying $A(s)\to 0$ as $s\to\infty$ and a second-order parameter $\rho\leq 0$, such that for all $x>0$
$$
\lim_{s\to \infty}\frac{1}{A(s)}
\left(
\frac{b(sx)-b(s)}{a(s)}-\int_1^{x}u^{\gamma-1}\diff u
\right)= \mathcal{H}_{\gamma,\rho}(x):=\int_1^{x}v^{\gamma-1}\int_{1}^v u^{\rho-1}\diff u \diff v;
$$
\item \label{en_sec_order_bias} The rate (or second-order auxiliary) function $A$ satisfies $\sqrt{k}A(m)\to\lambda\in \Real$.
\end{inparaenum}
\end{cond}
Conditions \ref{cond:second_order}\ref{en_sec_order_1}--\ref{en_sec_order_bias} are standard in the extreme values literature for bias quantification, and have been commonly used for decades in both cases of serial independence and dependence, see e.g. \citet{haan2006}, \citet{de2016adapting}, \cite{chavez2018extreme}, \citet{girard2021extreme}, \citet{davison2023tail}.

By \eqref{eq:DoA_TS}, we know that the distribution of the sample maximum—whether from a stationary or an i.i.d. sequence—when normalized by $a_{\theta}(m)$ and $b_{\theta}(m)$ and $a(\theta m)$ and $b(\theta m)$, respectively, converges to the same limiting distribution $G_\gamma$. Hence, to a first-order approximation, these pairs of normalizing sequences are asymptotically equivalent. However, at second order, a difference emerges in the residual errors, which can nonetheless still be quantified. The next result formalizes this argument.
\begin{lem} \label{lem:CONVNORMSEQ}
Under Condition \ref{cond:second_order} we have
$$
\sqrt{k}\left(\frac{a_{\theta}(m)}{a(\theta m)} - 1\right) = - \lambda  \frac{\theta^{\rho} - 1}{\rho} + o(1), \quad
$$
and
$$
\sqrt{k}\frac{b_{\theta}(m)-b(\theta m)}{a(\theta m)}= \lambda \theta^{\rho} \mathcal{H}_{\gamma, \rho} (\theta^{-1}) + \lambda  \frac{\theta^{\rho} - 1}{\rho} \frac{\theta^{-\gamma} - 1}{\gamma} + o(1).
$$
By the Slutsky’s lemma we then have
$$
\lim_{n\to\infty}
\Prob\left(M_{m}\leq a(\theta m)z +b(\theta m)\right)= G_{\gamma}(z).
$$
\end{lem}
In what follows, we assume for simplicity that the norming sequences are given by $a_m:=a(\theta m)$ and $b_m:=b(\theta m)$, so that the parameters $\mu$ and $\sigma$ correspond to these sequences.

Second, a key quantity in risk assessment is the return level associated with a return period $T$, often referred to as the $T$-time-horizon return level. Specifically, for a given time horizon $T=1,2,\ldots$, define the probability  $\tau\in(0,1)$ by $\tau=1-1/T$. The return level $R^{(m)}(\tau)$ is then the  $\tau$-quantile of the distribution of the sample maximum, obtained by inverting  $\tau=\Prob(M_m\leq x)$. In the case of independent observations, this quantity admits a natural interpretation. For instance, if $(X_1,\ldots,X_m)$ represent daily observations over a year (e.g., $m=366$), then $R^{(m)}(\tau)$ corresponds to the value expected to be exceeded by the annual maximum, on average, once every $T$ years. This interpretation does not hold in time series, however, thanks to \eqref{eq:DoA_TS} we still have like in the independence case (but with the obvious difference for the normalization constants)
\begin{equation}\label{eq:return_level}
R^{(m)}(\tau)\approx\mu + \sigma Q_\gamma(\tau),\quad m\to\infty,
\end{equation}
where the right-hand side above is the $\tau$-quantile of $G_\bfvartheta$ and $Q_\gamma(p)=G_\gamma^\leftarrow(p)$ for any $p\in[0,1]$, and where $Q_\gamma(p)=-\log(-\log(p))$ if $\gamma=0$. In the following sections, we develop the asymptotic theory for likelihood-based estimators of the return level.

Third, another important tool for risk assessment in time series is the VaR, a widely used risk measure in finance and economics. This measure, denoted by $Q(\tau)$, is indexed by a level $\tau \in(0,1)$, which specifies the frequency and severity of potential adverse events. From a probabilistic standpoint, VaR corresponds to the $\tau$-quantile of the stationary distribution $F$ of the time series.
A crucial challenge in risk protection is 
assessing the severity of rare events—those more extreme than any observed so far. This extrapolation task is inherently difficult because the true distribution $F$ is unknown in practice. However, EVT provides a framework to tackle this issue. Leveraging \eqref{eq:DoA_IID}–\eqref{eq:DoA_TS}, we can approximate
$\Prob(M_m\leq x)\approx F^{m\theta}(x)\approx G_\bfvartheta(x)$ as $m\to\infty$. Therefore, using the identity $\tau=F(x)$ and applying the EVT-based approximation, the VaR at an extreme level (i.e., as $\tau\to 1$) can be approximated by
\begin{equation}\label{eq:ext_quantile}
Q(\tau)\approx \mu + \sigma \,Q_\gamma(\tau^{m\theta}).
\end{equation}

We complete this section introducing the following conditions that specify the type of serial dependence structure required for the class of time series we consider.
\begin{cond}[Mixing conditions]\label{cond:ts_cond}
The stationary time series $(X_t)_{t\geq 1}$ satisfies:
\begin{inparaenum}
\item \label{cond_alpha_mix1} For any $\varepsilon\in(0,1]$, Let $\mathcal{F}_{i,j}^{\varepsilon}$ be the sigma-field generated by  $(F(X_i) \ind (F (X_i) > 1 - \varepsilon),\ldots,F(X_j) \ind (F (X_j) > 1 - \varepsilon))$. Let $\alpha_\varepsilon:\Nat\to[0,1]$ be 
the mixing coefficients defined by $\alpha_\varepsilon(l)=\sup_{r\in\Nat} \alpha_r^\varepsilon (l)$, where
$$
\alpha_r^\varepsilon(l)=\sup_{A\in\mathcal{F}_{1,r}^{\varepsilon}, B\in\mathcal{F}_{r+l,\infty}^{\varepsilon}}|\Prob(A \cap B)-\Prob(A)\Prob(B)|.
$$
Then, we have $\alpha_\varepsilon(l)\to 0$ as $l\to\infty$, which is a weaker assumption than the standard alpha-mixing condition (recovered when $\varepsilon=1$).
\item \label{cond_alpha_mix2} The coefficients in \ref{cond_alpha_mix1} also satisfy 
$$
\lim_{n\to\infty} k^2 (\log k)^4\sum_{r=1}^{k} r \alpha_\varepsilon((r-1)(m+l)+l)=0,
$$
where, for $n=1,2,\ldots$, $m=m_n$, $l=l_n$ and $k=k_n$ are sequences such that $l\to\infty$, $m\to\infty$ and $k\to\infty$ as $n\to\infty$ and $m=o(n)$, $k=o(n)$ and $l=o(m)$. The sequences $(m_n)_{n\geq 1}$, $(l_n)_{n\geq 1}$ and $(k_n)_{n\geq 1}$ are known as big-block, small-block and intermediate sequences, respectively. 
\item\label{cond_beta_mix} Let $(\beta_n)_{n\geq 1}$ be the sequence defined by
\begin{equation}\label{eq:beta_coef}
\beta_n = \sup_{s\in[0,1]} \biggr\vert \Prob \left(\max_{1 \leq i \leq m} F^{m\theta}(X_i)\leq s \right) - \Prob\left(\max_{1\leq i \leq \lfloor \theta m \rfloor} F^{\lfloor \theta m \rfloor}(\widetilde{X}_i)\leq s\right) \biggr\vert.
\end{equation}
For simplicity set $\beta=\beta_n$. Then, $k\beta\to 0$ as $n\to\infty$.
\end{inparaenum}
\end{cond}
Assumption \ref{cond:ts_cond}\ref{cond_alpha_mix1} is weaker than the usual alpha-mixing condition \citep{berghaus2018}, with the latter that is a standard argument in the literature on mixing stationary series and is widely adopted in the extreme values literature (e.g. \citealp{leadbetter1983, leadbetter1988extremal, rootzen2009weak}). 
Assumption \ref{cond:ts_cond}\ref{cond_alpha_mix2} ensures that the dependence between two sample maxima computed over a block of $m$ variables, that are separated by a sequence of $l$ variables, becomes asymptotically negligible. 
Assumption \ref{cond:ts_cond}\ref{cond_beta_mix} guarantees that a sample maximum computed over a block of $m$ variables, behaves asymptotically (in the sense of \eqref{eq:beta_coef}) as a sample maximum computed over a block of independent variables of size $\lfloor \theta m \rfloor$.
\begin{rem} \label{rem:betamixing}
If the time series satisfies an analogue of Condition \ref{cond:ts_cond}\ref{cond_alpha_mix1} with mixing coefficients $\phi_\varepsilon:\Nat\to[0,1]$, given by $\phi_\varepsilon(l)=\sup_{r\in\Nat} \phi_r^\varepsilon (l)$, for any $\varepsilon>0$, where
$$
\phi_r^\varepsilon(l)=\Expect(\sup(|\Prob(B\mid \mathcal{F}_{1,r})-\Prob(B)|: B\in \mathcal{F}_{r+l,\infty})),
$$
which is the absolutely regular or beta-mixing condition when $\varepsilon=1$, then the Condition \ref{cond:ts_cond}\ref{cond_alpha_mix2} simplifies to $\lim_{n \rightarrow \infty} k\phi_\varepsilon (l)=0$, 
see Remark A.11 of the Supplement for details.
\end{rem}

%
\section{Likelihood asymptotic theory} \label{sec:LAT}
%
%
From a modeling view point, the loglikelihood function corresponding to one observation of the GEV model $G_{\bfvartheta}$, with parameters $\bfvartheta=(\gamma, \mu, \sigma)^\top$, is
$$
\ell_{\bfvartheta}(x)=
\begin{cases}
	-\log \sigma -\log g_\bfvartheta(x),& \text{if } 1+\frac{\gamma}{\sigma}x>0\\
	-\infty,& \text{otherwise}.
\end{cases}
$$
and the related {\it expected information} function is
$$
I(\bfvartheta)=-\int_0^1\frac{\partial^2 \ell_{\bfvartheta}}{\partial \bfvartheta\partial \bfvartheta^\top}\left(\mu + \sigma\frac{(-\log u )^{-\gamma}-1}{\gamma}\right)\diff u, \quad \forall\, \bfvartheta\in\Real^2\times(0,\infty).
$$
For a given $\bfvartheta_0$, the expected information $\bfI(\bfvartheta_0)$ (see Supplement, Section D) is positive definite for $\gamma_0>-1/2$ and so in the sequel we work with $\Theta=(-1/2,\infty)\times\Real\times(0,\infty)$.

Now, let $\bfX_n=(X_1,\ldots,X_n)$ be a sample from the stationary time series $(X_t)_{t\geq 1}$, with marginal distribution $F_0$, satisfying $F_0\in\mathcal{D}(G_{\gamma_0})$, with $\gamma_0>-1/2$. The theory of Section \ref{eq:backg} gives rise to one of the most popular inferential methods in the univariate extreme values context known as {\it Block Maxima} (BM). According to the argument on ``big-blocks of size $m$ separated by small blocks of size $l$'' the BM method suggests to split $\bfX_n$ into $k$ blocks of $m+l$ consecutive variables, where $k=\lfloor n/(m+l)\rfloor$, and then compute the series of disjoint-block sample maxima 
\begin{equation} \label{eq:defbm}
	M_{i,m}=\max(X_{(1+(i-1)(m+l))},\ldots,X_{m+(i-1)(m+l)}), \quad i=1,\ldots,k.
\end{equation}
We assume that the integers $k$, $l$, $m$ and $n$ behaves according to Condition  \ref{cond:ts_cond}\ref{cond_alpha_mix2}. 
By Lemma \ref{lem:CONVNORMSEQ}, for sufficiently large $m$, the sequence $M_{m,i}$ is approximately distributed according to the GEV model $G_{\bfvartheta_0}$, where $\bfvartheta_0=(\gamma_0,\mu_0, \sigma_0)^\top$, and where $\mu_0$ and $\sigma_0$ are representative of the norming constants $b_{m,0}$ and $a_{m,0}$, respectively. The independence case is recovered setting $l=0$. 
\begin{rem} \label{rem:smallblocks}
	If the small block size $l$ grows slowly enough so that the first $m$ variables of a block of size $m+l$ generates the largest observation  with probability tending to one, then our theory could equivalently be reformulated by computing block maxima over blocks of size $m+l$.
	A sufficient condition (see Lemma A.24) for this to happen is
	\begin{equation*}
		\lim_{n \rightarrow \infty} \Prob \left( \bigcap_{i = 1}^{k} \left\{ M_{i, m} \geq \max_{i (m + l) - l + 1\leq j \leq i (m + l)} X_j \right\}  \right) = 1.
	\end{equation*}
\end{rem}

We now present our first main theoretical result, extending the tail quantile process' asymptotic representation of \citet{ferreira2015} from the independence case to the time series setting, which will later serve as the foundation for many of our key inferential developments.
Let $M_{1,k;m}\leq\cdots\leq M_{k,k;m}$ be the $k$ order statistics of the sequence of block maxima. Let $(\mathbb{Q}_{k,m}(s), 0<s<1)$ be the quantile process defined by
$$
\mathbb{Q}_{k,m}(s)=\sqrt{k}\left(\frac{M_{\lceil sk \rceil:k,m}-b_{m,0}}{a_{m,0}}-Q_{\gamma_0}(s)\right).
$$
\begin{theo}[Tail quantile process] \label{theo:TQP}
	Under Conditions  \ref{cond:second_order} and \ref{cond:ts_cond} there exists a sequence of Brownian bridges $(B_n)_{n\geq 1}$ such that for any $0<\varepsilon<1/2$,
	$$
	\mathbb{Q}_{k,m}(s)=\frac{B_n(s)}{s(-\log s)^{1+\gamma_0}}+\sqrt{k}A(\lfloor \theta_0 m \rfloor)\mathcal{H}_{\gamma_0,\rho_{0}}\left(\frac{1}{-\log s}\right)+ R_n(s),
	$$
	where
	$$
	R_n(s)=(s^{-1/2-\varepsilon} (1-s)^{-1/2-\gamma_0-\rho_{0}-\varepsilon})o_{\Prob}(1),
	$$
	$\sqrt{k}A(\lfloor \theta_0 m \rfloor) \rightarrow \lambda_0 \theta_{0}^{\rho_{0}}$ as $n \rightarrow \infty$,
	and $o_{\Prob}(1)$ is uniform for $1/(1+k)\leq s \leq k/(1+k)$.
\end{theo}
\begin{rem}\label{rem:weak_SOC}
	The expansion in Theorem \ref{theo:TQP}  is also valid with $\mathcal{H}_{\gamma_0,\rho_0} \equiv 0$ and a suitable choice of the scaling sequence $a(\cdot)$ 
	under an alternative assumption to the classical second order condition. Precisely, the latter can be replaced by
	\begin{equation}\label{eq:weaker}
		\sup_{x \geq 1/(1+\epsilon) } \sqrt{k_0}  \left|
		\frac{b(m_0 x)-b(m_0)}{a(m_0)} - \int_1^x u^{\gamma_0-1} \diff u
		\right|x^{-\gamma_0-1/2} =o(1),
	\end{equation}
	for some $\epsilon>0$ and with $m_0 = m/ \log k_0$, $k_0 = -k /W(-1/k)$, where $W$ is the product log function; see Section B of the Supplement for an extended formal statement and proof. Unlike Condtion \ref{cond:second_order}\ref{en_sec_order_1},  there always exist some sequences $m =o(n)$ and $k\sim n/m$ complying with \eqref{eq:weaker}. 
	In simple terms, as long as the block size is sufficiently large (and Condition \ref{cond:ts_cond} is fulfilled), it is possible to approximate the tail quantile process by a weighted Brownian bridge even for those time series whose stationary distribution is not available in closed form, and for which directly verifying Condtion \ref{cond:second_order}\ref{en_sec_order_1} is more challenging.
	Important examples of such time series are solutions of stochastic recurrence equations, as e.g. squared ARCH processes.
	See \cite{drees2000, drees2003} for similar considerations on the quantile process of threshold exceedances. 
\end{rem}

Now we develop the asymptotic theory the likelihood-based inference. Let 
\begin{equation}\label{eq:logllik}
	\mathscr{L}_n(\bfvartheta)=\frac{1}{k}\sum_{i=1}^k \ell_{\bfvartheta}(M_{i,m})
\end{equation}
be the empirical pseudo loglikelihood process relative to $\bfvartheta\in\Theta$.
Note that the supremum of loglikelihood process $\mathscr{L}_n(\bfvartheta)$ over $\bfvartheta \in\Theta$ is unbounded for $\gamma$ that goes to infinity, see Section A.3.3 of the supplement for details. To remedy this we need to make $\Theta$ sample size dependent in $\gamma$, namely we restrict the parameter space to $\Theta_{n} = \{\bfvartheta \in \Theta : \gamma < k^{q} \}$ for any $q < 1$. Now, we define the maximum likelihood estimator (MLE) by
\begin{equation}\label{eqMLE}
	\widehat{\bfvartheta}_n\in
	\argmax_{\bfvartheta\in\Theta_{n}} \mathscr{L}_n(\bfvartheta),
\end{equation}
where we omit the dependence on $q$ to ease notation.
Asymptotic theory requires we let $n$ go to infinity, allowing $m$ to vary accordingly (see Condition \ref{cond:ts_cond}\ref{cond_alpha_mix2}). Consequently, the parameters $\mu$ and $\sigma$, representing the norming sequences $b_{m,0}$ and $a_{m,0}$, also vary, which complicates the resulting estimation theory. To work out the theory we consider the reparametrization $\bftheta=r(\bfvartheta)=(\gamma,(\mu-b_{m,0})/a_{m,0}, \sigma/a_{m,0})^\top$ for all $\bfvartheta\in\Theta_n$ and the alternative empirical loglikelihood $L_n(\bftheta)=k^{-1}\sum_{i=1}^k \ell_{\bftheta}(\widetilde{M}_{i,m})$, based on the so-called pseudo-``observations'' $\widetilde{M}_{i, m}=(M_{i,m}-b_{m,0})/a_{m,0}$ for $i=1,\ldots,k$ \citep[see e.g.][for the independence case]{dombry2019, dombry2023}.  
We clarify to avoid confusion about the notation for the remainder of the article that $\theta$ is the extremal index, while $\bftheta$ denotes instead the vector of the GEV standardized parameters. 
Once the asymptotic behaviour for $\widehat{\bftheta}_n$ is established, that for $\widehat{\bfvartheta}_n=r^{-1} (\widehat{\bftheta}_n)$ is obtained as byproduct. Note that $\bftheta_0=(\gamma_0,0,1)^\top$ and $L_n(\bftheta)=\mathscr{L}_n(\bfvartheta)+\log(a_{m,0})$.

Before stating our main likelihood results we specify some useful notation.
We recall that, for every $\bftheta\in\Theta_n$, $\bfS_n(\bftheta)=(\partial /\partial \bftheta) L_n(\bftheta)$ and $\bfJ_n(\bftheta)=(\partial^2/(\partial \bftheta \partial \bftheta^\top)) L_n(\bftheta)$ are the rescaled {\it score} and rescaled negative {\it observed information} functions pertaining to $L_n(\bftheta)$, respectively. For  $\bftheta_0\in\Theta_n$ let $\mathbb{B}_{\varepsilon}(\bftheta_0)=(\bftheta\in\Theta: \|\bftheta-\bftheta_0\|< \varepsilon)$ be the ball with center $\bftheta_0$ and radius $\varepsilon>0$.
 Finally, $ \mathcal{N}(\bfm,\bfSigma)$ denotes a multivariate normal distribution with mean vector $\bfmu$ and covariance matrix $\bfSigma$, which reduces to $N(\mu,\sigma^2)$ in the univariate case. We denote by $ \mathcal{N}(\cdot;\bfmu,\bfSigma)$ the related probability measure.
\begin{theo}[Likelihood expansions]\label{theo:loglike_expansion}
	Under Conditions \ref{cond:second_order} and \ref{cond:ts_cond} the following results hold.
	\begin{itemize}
		\item Let $C=C_n$ with $C\to\infty$ as $n\to\infty$ and $C=O(k^\delta)$, with $\delta < \min(0.5, \gamma_0 + 0.5)$. For $\varepsilon= C/\sqrt{k}$ we have
		$$
		\sup_{\bftheta\in\mathbb{B}_{\varepsilon}(\bftheta_0)}\|\bfJ_n(\bftheta)+\bfI(\bftheta_0)\|=o_{\Prob}(1),
		$$
		\item For all $c>0$ and $\varepsilon= c/\sqrt{k}$ we have
		$$
		\sup_{\bftheta\in\mathbb{B}_{\varepsilon}(\bftheta_0)}
		\left|
		L_n(\bftheta)-L_n(\bftheta_0)-(\bftheta-\bftheta_0)^\top \bfS_n(\bftheta_0)
		+\frac{(\bftheta-\bftheta_0)^\top \bfI(\bftheta_0) (\bftheta-\bftheta_0)}{2} 
		\right|=o_{\Prob}\left(\frac{1}{k}\right).
		$$
		and
		$$
		\sup_{\bftheta\in\mathbb{B}_{\varepsilon}(\bftheta_0)}
		\left|
		\bfS_n(\bftheta)- \bfS_n(\bftheta_0) - \bfI(\bftheta_0) (\bftheta-\bftheta_0)
		\right|=o_{\Prob}\left(\frac{1}{\sqrt{k}}\right).
		$$ 
		\item We have
		$$\sqrt{k}\bfS_n(\bftheta_0)\stackrel{d}{\rightarrow} \mathcal{N}_3(\lambda_0\theta_0^{\rho_{0}}\bfb_0, \bfI(\bftheta_0)),$$
		where $\bfb_0$ is given in Section D of the Supplement.
	\end{itemize}
\end{theo}
\begin{theo}[Local and global bounds]\label{theo:local_global_bounds}
	Under Conditions \ref{cond:second_order} and \ref{cond:ts_cond} there exist $\varepsilon>0$ and constants $c_1, c_2 >0$ such that with probability tending to one  as $n\to\infty$
	$$
	L_n(\bftheta)-L_n(\bftheta_0) \leq (\bftheta-\bftheta_0)^\top \bfS_n(\bftheta_0)-c_1 \|\bftheta-\bftheta_0\|^{2}, \quad \forall\, \bftheta\in\mathbb{B}_{\varepsilon}(\bftheta_0)
	$$
	and 
	\begin{eqnarray*}
		L_n(\bftheta)-L_n(\bftheta_0) &\leq& -c_2, \quad \forall \,\bftheta\in\Theta_n\backslash \mathbb{B}_{\varepsilon}(\bftheta_0).\\
	\end{eqnarray*}
\end{theo}
\begin{cor}[Uniqueness and asymptotic normality of the MLE]\label{cor:mle}
	Under Conditions \ref{cond:second_order} and \ref{cond:ts_cond} the following results hold.
	\begin{itemize}
		\item The MLE in \eqref{eqMLE} is the unique global maximiser of the loglikelihood function $\mathscr{L}_n({\bfvartheta})$ over $\Theta_n$ with probability tending to one, as $n\to\infty$.
		\item For $\bfb_0$ as in Theorem \ref{theo:loglike_expansion},
		$$
		\sqrt{k}\left(
		\widehat{\gamma}_n-\gamma_0, \frac{\widehat{\mu}_n- b_{m,0}}{a_{m,0}}, \frac{\widehat{\sigma}_n}{a_{m,0}}
		-1\right)\stackrel{d}{\rightarrow} \mathcal{N}(\lambda_0\theta_0^{\rho_{0}}\bfI^{-1}(\bftheta_0)\bfb_0, \bfI^{-1}(\bftheta_0)).
		$$
	\end{itemize}
\end{cor}
Based on Corollary \ref{cor:mle}, an approximation of an asymptotic $(1-\alpha)$-equi-tailed confidence interval estimator for $\gamma_0$  is given by
\begin{equation}\label{eq:CI_MLE_gamma}
	\left[
	\widehat{\gamma}_n-\widehat{b}_n\pm \frac{z_{1-\alpha/2}\widehat{\psi}_n}{\sqrt{k}}
	\right],
\end{equation}
where $z_{1-\alpha/2}$ is the $(1-\alpha/2)$-quantile of the standard normal distribution. In particular,  $\widehat{b}_n$ is an estimator of $b_0$, i.e. the first element of $k^{-1/2}\lambda_0 \theta_{0}^{\rho_0} \bfI^{-1}(\bftheta_0) \bfb_0$. $\widehat{\psi}_n$ is the square root of the first diagonal element of $\widehat{\bfI}_{n}^{-1} = (-\bfA_n  \mathcal{\bfJ}_n(\widehat{\bfvartheta}_{n}) \bfA_n)^{-1}$ where  $\mathcal{\bfJ}_n(\widehat{\bfvartheta}_n)=(\partial^2/(\partial \bfvartheta \partial \bfvartheta^\top)) \mathscr{L}_n(\bfvartheta)|_{\bfvartheta=\widehat{\bfvartheta}_n}$ and $\bfA_n$ is a diagonal matrix with elements $(1, \widehat{\sigma}_n, \widehat{\sigma}_n)^{\top}$, which is an estimator of $\psi_0$, i.e. the square root of the first diagonal element of $\bfI^{-1}(\bftheta_0)$.

Risk assessment based on $T$-return level $R^{(m)}(\tau)$, for $\tau\in(0,1)$ is widely used in practice, particularly under the assumption of independence due to its intuitive interpretation. In the time series context, EVT builds asymptotic models for sample maxima by leveraging the fact that maxima computed over big-blocks of size $m$, separated by small-blocks of size $l$, become asymptotically independent as both $m$ and $l$ grow with the overall sample size $n$ (e.g., \citealp[][Ch. 3]{leadbetter1983}, \citealp[][Ch. 10]{beirlant2006}). As a result, for sufficiently large samples, the sample maxima can be regarded  as (approximately) independent, which supports retaining the conventional (in this case heuristic) interpretation of $R^{(m)}(\tau)$ even in time series \citep[e.g.][p. 381]{beirlant2006}. From a practical standpoint, it may be 
convenient to fit $G_\bfvartheta$ to $M_{i,m}$, $i=1,\ldots,k$, i.e. sample maxima of block size $m$ (e.g., annual maxima), however, for interpretational purposes it may be preferable to assess the return level $R^{(m^\star)}(\tau)$ relative to $M_{i,m^\star}$, whose block size is $m^\star\geq m$ (e.g., decennial maxima). In this case, leveraging formula \eqref{eq:return_level}, for large $m$, the $T$-return level becomes $R^{(m^\star)}_{0}(\tau)\approx R^{(m)}_{0}(\tau^{m / m^\star})\approx \mu_0 +\sigma_0 Q_{\gamma_0}(\tau^{m / m^\star})$ and its MLE is $\widehat{R}_n(\tau^{m / m^\star})=\widehat{\mu}_n+\widehat{\sigma}_n Q_{\widehat{\gamma}_n}(\tau^{m/m^\star})$.
\begin{theo}[Return Levels] \label{thm:RETURNLV}
	Work under the Conditions of Theorem \ref{theo:local_global_bounds}. 
	Let $\tau \in (0, 1)$ and $m^\star = m^\star_n$ such that $m^\star / m \rightarrow c\in[1,\infty)$ as $n \rightarrow \infty$. Then,  $m^\star/(-m\log \tau) \to \omega\in(0,\infty)$ as $n \rightarrow \infty$ and
	\begin{equation*}
		\frac{\sqrt{k}(\widehat{R}_n(\tau^{m/m^\star}) - R^{(m^\star)}_{0}(\tau))}{a_{m,0} q_{\gamma_0} (\tau^{m/m^\star})} 
		\overset{d}{\rightarrow} 
		\mathcal{N}\left(\lambda_0\theta_{0}^{\rho_{0}} \left(\mathbf{v}_{\gamma_0}^{\top} \bfI^{-1}(\bftheta_0) \bfb_0 + b_{\gamma_0}\right), \mathbf{v}_{\gamma_0}^{\top} \bfI^{-1}(\bftheta_0) \mathbf{v}_{\gamma_0} \right),
	\end{equation*}
	where $q_{\gamma_0}(x)=\partial/(\partial \gamma) Q_{\gamma}(x)$, $\bfb_0$ is as in Theorem \ref{theo:loglike_expansion}, 
	\begin{equation*}
		\mathbf{v}_{\gamma_0} = \left( 
		1, 
		\frac{1}{q_{\gamma_0} (\tau^{m/m^\star})}, 
		\frac{Q_{\gamma_0} (\tau^{m/m^\star})}{q_{\gamma_0} (\tau^{m/m^\star})} \right)^\top, 
		\quad
		b_{\gamma_0} =  -\frac{\mathcal{H}_{\gamma_0, \rho_0}(\omega)+\omega^{\gamma_0 + \rho_0}}{q_{\gamma_0} (\tau^{m/m^\star})}.
	\end{equation*}
	If instead $m^{\star} / m \rightarrow \infty$ as $n \rightarrow \infty$, then the result above still holds with $\mathbf{v}_{\gamma_0} =(1, 0, 0)^\top$ if $\gamma_0 \geq 0$ and  $\mathbf{v}_{\gamma_0}=(1, \gamma_{0}^{2}, -\gamma_0)^{\top}$ otherwise, $b_{\gamma_0}=-\gamma_0/(\gamma_0+\rho_0)$ if $\gamma_0<0$ and $b_{\gamma_0}=0$ otherwise, provided that $\log (m^\star/(-m\log \tau)) = o(\sqrt{k})$ and $\rho_{0} < 0$ in the case $\gamma_0 \geq 0$.
\end{theo}
Leveraging Theorem \ref{thm:RETURNLV}, an approximation of an asymptotic $(1-\alpha)$-equi-tailed confidence interval estimator for $R^{(m^\star)}_{0}(\tau)$ is given by
\begin{equation}\label{eq:CI_MLE_RL}
	\left[
	\widehat{R}_n(\tau^{m / m^\star})-\widehat{\sigma}_n q_{\widehat{\gamma}_n}(\tau^{m / m^\star})\left(
	\widehat{b}_n\pm \frac{z_{1-\alpha/2}\widehat{\psi}_n}{\sqrt{k}}
	\right)
	\right],
\end{equation}
where $\widehat{b}_n$ is an estimator of $b_0=k^{-1/2}\lambda_0\theta_{0}^{\rho_{0}} \left(\mathbf{v}_{\gamma_0}^{\top} \bfI^{-1}(\bftheta_0) \bfb_0 + b_{\gamma_0}\right)$, $\widehat{\sigma}_n$ is the MLE of $\sigma_0$,  $q_{\widehat{\gamma}_n}$ is the estimator of $q_{\gamma_0}$ obtained pugging-in the MLE $\widehat{\gamma}_n$ of  $\gamma_0$, and  
 $\widehat{\psi}_n = (\mathbf{v}_{\widehat{\gamma}_n}^{\top} \widehat{\bfI}_{n}^{-1} \mathbf{v}_{\widehat{\gamma}_n})^{1/2}$ is an estimator of $\psi_0=(\mathbf{v}_{\gamma_0}^{\top} \bfI^{-1}(\bftheta_0) \mathbf{v}_{\gamma_0})^{1/2}$.

Risk assessment for events more extreme than those yet seen can be carried out using the VaR at an extreme level $\tau$ close to 1, effectively extrapolating into the far tail of the stationary distribution of the time series. Let $\tau_E=\tau_{E,n}$ denote such an extreme level, where $\tau_E\to1$ as $n\to\infty$, and this convergence is sufficiently fast in the sense that $m(1-\tau_E)=o(1)$, e.g. $\tau_E=1-c/n$ for some $c>0$. Under these conditions, and following equation \eqref{eq:ext_quantile}, the VaR at level $\tau_E$ can be approximated by $Q_0(\tau_E)\approx \mu_0 +\sigma_0 Q_{\gamma_0}(\tau_E^{m\theta_0})$ with its MLE given by $\widehat{Q}_n(\tau_E)=\widehat{\mu}_n+\widehat{\sigma}_n Q_{\widehat{\gamma}_n}(\tau_E^{m\widehat{\theta}_n})$, where $\widehat{\theta}_n$ is an estimator of $\theta_0$, \citep[e.g.,][]{robert2009, berghaus2018}.
\begin{theo}[Extreme VaR] \label{thm:EXTRQUANT}
	Let  $\tau_E \to 1$ as $n\to\infty$ and such that  $m(1-\tau_E)=o(1)$.
	Work under the Conditions of Theorem \ref{theo:local_global_bounds} and suppose that both $-\log(-\theta_0 m \log \tau_E)=o(\sqrt{k})$ and $\rho_0<0$ if $\gamma_0 \geq 0$.
	Let $\widehat{\theta}_n$ be a $\sqrt{k}$-consistent estimator of $\theta_0$. Then, 
	\begin{equation*}
		\frac{ \sqrt{k} (\widehat{Q}_n(\tau_E) - Q_0(\tau_E))}{a_{m,0} q_{\gamma_0} (\tau_E^{m\theta_0})} 
		\overset{d}{\rightarrow} 
		\mathcal{N} \left(\lambda_0\theta_{0}^{\rho_{0}} \left(\mathbf{v}_{\gamma_0}^{\top} \bfI^{-1}(\bftheta_0) \bfb_0 + b_{\gamma_0}\right), \mathbf{v}_{\gamma_0}^{\top} \bfI^{-1}(\bftheta_0) \mathbf{v}_{\gamma_0} \right),
	\end{equation*}
	where $q_{\gamma_0}$ and $\bfb_0$ are as in Theorem \ref{thm:RETURNLV}, $\mathbf{v}_{\gamma_0} =(1, 0, 0)^\top$ if $\gamma_0 \geq 0$ and   $\mathbf{v}_{\gamma_0}=(1, \gamma_{0}^{2}, -\gamma_0)^{\top}$ otherwise, $b_{\gamma_0}=-\gamma_0/(\gamma_0+\rho_0)$ if $\gamma_0<0$ and $b_{\gamma_0}=0$ otherwise.
\end{theo}

From a practical view point, to achieve good performance with finite samples,  the derivation of confidence intervals for $Q_0(\tau_E)$ requires an estimator $\widehat{\theta}_n$ of $\theta_0$ and of its variance. Our choice for $\widehat{\theta}_n$ is the MLE with disjoint-block maxima discussed in \citet[][]{northrop2015efficient, berghaus2018}, that we briefly summarize here. By formula \eqref{eq:defbm} we define the sample of maxima $M_{i,\lm}$, $i=1,\ldots,\lk$,  where $l = 0$, $\lm$ is not necessarily equal to $m$ but still $\lm=\lm_n$ with $\lm=o(n)$. From \eqref{eq:DoA_extremal_index} we have that for any $x\in(0,1)$, $\lm\overline{F}(u_{\lm})\to x$  and $\Prob(M_{\lm}\leq u_{\lm})\to e^{-\theta_0 x}$ as $\lm\to\infty$, as soon $u_{\lm}=Q(e^{-x/\lm})$.
Since $F$ is a nondecreasing function, then $F(M_{\lm})=\max(U_{1},\ldots, U_{\lm})$, where $U_i=F(X_i)$ is standard uniform random variable. 
Then, the sequence of random variables $Y_{i,\lm}=-\lm\log(F(M_{i,\lm}))$, $i=1,\ldots,\lk$, are approximately exponentially distributed with mean $1/\theta_0$, as $n\to\infty$, since that $\Prob(Y_{i,\lm}\geq x)=\Prob(M_{i,
	\lm}\leq u_{\lm})$. Because $F$ is unknown, then $Y_{i,\lm}$ can be approximated by $\widehat{Y}_{i,\lm}=-\lm\log(F_n(M_{i,\lm}))$, where $F_n(x)=n^{-1}\sum_{t=1}^n\ind(X_t\leq x)$ 
and $\ind(E)$ is that is the indicator function of the event $E$. We define the pseudo empirical log-likelihood process
\begin{equation}\label{eq:ext_ind_likelihood}
	\mathscr{L}_n(\theta)=\left(\log \theta - \frac{\theta}{\lk}\sum_{i=1}^{\lk} \widehat{Y}_{i,\lm}\right)
\end{equation}
and the corresponding MLE 
\begin{equation}\label{eq:ext_ind_MLE}
	\widehat{\theta}_n= \argmax_{\theta\in(0,1]} \mathscr{L}_n(\theta).
\end{equation}
For the estimation of the asymptotic variance of $\widehat{\theta}_n$,  which, according to \citet{berghaus2018}, is  $\theta_0^4\sigma^2(\theta_0)$, we propose to estimate $\sigma^2(\theta_0)$ via 
\begin{equation}\label{eq:est_var_est_theta}
	\widetilde{\sigma}_{n}^{2}=\frac{1}{K}\left( \widetilde{\sigma}_{n, 1}^{2} +  \sum_{i = 2}^{K} \widetilde{\sigma}_{n, \lceil \frac{i - 1}{K} \lm \rceil + 1 }^{2} \right), \quad K \geq 1,
\end{equation}
with
\begin{align*}
	\widetilde{\sigma}_{n, 1}^2   &= \frac{1}{\lk}  \sum_{i = 1}^{\lk} \Biggr( \widehat{Y}_{i, \lm} - \frac{1}{\lk} \sum_{l = 1}^{\lk}  \widehat{Y}_{l, \lm} 
	+ 
	\sum_{1 + (i - 1) \lm \leq s \leq i \lm} \frac{1}{\lk} \sum_{l = 1}^{\lk} \left\{ \frac{F_{n} (M_{l, \widetilde{m}}) - \ind \left( F_{n} (X_s) \leq F_{n} (M_{l, \widetilde{m}}) \right)}{F_{n} (M_{l, \widetilde{m}})} \right\} \Biggr)^2\\
	\widetilde{\sigma}_{n, j}^2  &= \frac{1}{\lk - 1}  \sum_{i = 1}^{\lk - 1} \Biggr(\widehat{Y}_{i, \lm}^{(j)} - \frac{1}{\lk - 1} \sum_{l = 1}^{\lk - 1}  \widehat{Y}_{l, \lm}^{(j)}\\ 
	& + 
	\sum_{j + (i - 1) \lm \leq s \leq j - 1 +  i \lm} \frac{1}{\lk - 1} \sum_{l = 1}^{\lk -1} \Biggr\{ \frac{F_{n}^{(j)} (M_{l, \widetilde{m}}^{(j)}) - 
		\ind \left( F_{n}^{(j)} (X_{s}) \leq F_{n}^{(j)} (M_{l, \widetilde{m}}^{(j)}) \right)}{F_{n}^{(j)} (M_{l, \widetilde{m}}^{(j)})} \Biggr\} \Biggr)^2
\end{align*}
and where $F_{n}^{(j)} (x) = (n - j + 1)^{-1} \sum_{t = j}^{n} \ind (X_t \leq x)$ and $M_{l, \lm}^{(j)}$ are the empirical distribution function and the $l$th sample maximum of block size $\lm$ of $X_j, \dots, X_n$,  respectively, and $\widehat{Y}_{i, \lm}^{(j)} = - \lm (\log F_{n}^{(j)} (M_{i, \lm}^{(j)}))$. 
This estimator is an alternative version of the estimator proposed by \citet{berghaus2018} that shows to be more efficient, see Section C of the Supplement for the details.
Now, exploiting Theorem \ref{thm:EXTRQUANT}, an approximation of an asymptotic $(1-\alpha)$-equi-tailed confidence interval estimator for $Q_{0}(\tau_E)$ is given by
\begin{equation}\label{eq:CI_MLE_EXTQ}
	\left[
	\widehat{Q}_n(\tau_E)-\widehat{\sigma}_n q_{\widehat{\gamma}_n}(\tau^{m \widehat{\theta}_n})\left(
	\widehat{b}_n\pm \frac{z_{1-\alpha/2}\widehat{\psi}_n}{\sqrt{k}}
	\right)
	\right],
\end{equation}
where $\widehat{b}_n$ is an estimator of $b_0=k^{-1/2}\lambda_0\theta_{0}^{\rho_{0}} \left(\mathbf{v}_{\gamma_0}^{\top} \bfI^{-1}(\bftheta_0) \bfb_0 + b_{\gamma_0}\right)$ and $\widehat{\psi}_n = (\mathbf{v}_{\widehat{\gamma}_n}^{\top} \widehat{\bfI}_{n}^{-1} \mathbf{v}_{\widehat{\gamma}_n} + \widehat{\varsigma}_n)^{1/2}$ is an estimator $\psi_0=(\mathbf{v}_{\gamma_0}^{\top} \bfI^{-1}(\bftheta_0) \mathbf{v}_{\gamma_0})^{1/2}$, and where the correction term 
$$
\widehat{\varsigma}_n=\frac{\widetilde{\sigma}_{n}^{2} \widehat{\theta}_{n}^2 }
{\left(\left(\log (\tau_{E}^{m \widehat{\theta}_n}))^{\widehat{\gamma}_n} q_{\widehat{\gamma}_n} (\tau_{E}^{m \widehat{\theta}_n}\right)\right)^{2}}
$$ 
helps to obtain better coverage with finite-samples, while it vanishes when $n\to\infty$.

The coverage probability of the confidence interval estimators \eqref{eq:CI_MLE_RL} and \eqref{eq:CI_MLE_EXTQ} can asymptotically achieve the nominal level, even neglecting the estimator $\widehat{b}_n$ when the bias term $\lambda_0=0$. However, for finite samples of small or moderate size, the performance of such intervals can be poor when the distribution of the MLE deviates from normality. This limitation arises because the results in Theorems \ref{thm:EXTRQUANT} and \ref{thm:RETURNLV} rely not only on the tools developed in Corollary \ref{cor:mle}, but also on a first-order delta method approximation. In alternative, we propose to use the following asymmetric confidence intervals that do not rely on the delta method and demonstrate that, under the condition $\lambda_0=0$, their coverage probability also converges asymptotically to the nominal level. 

\begin{cor}[Asymmetric confidence intervals] \label{cor:CONFINT}
	For any measurable set $B_1\subset \Theta_n$, define the sequence of random probability measures 
	\begin{equation*}\label{eq:Gaussian_measure}
		\varGamma_{n} (B_1) 
		= 
		\mathcal{N} \left(B_1; \widehat{\bfvartheta}_n, -k^{-1}\mathcal{\bfJ}_n(\widehat{\bfvartheta}_n)\right).
\end{equation*}
	\begin{itemize}
		\item Work under the conditions of Theorem \ref{thm:EXTRQUANT} with $\lambda_0 = 0$. For any measurable set $B_2\subset (0,1]$, define the sequence of random probability measures 
		\begin{equation*}\label{eq:extremal_index_measure}
			\varLambda_n(B_2)
			= 
			\mathcal{N} (B_2; \widehat{\theta}_n, k^{-1} V_n),
		\end{equation*}
		where $V_n=v+o_{\Prob}(1)$, for some $v>0$. 
		For every measurable set $B\subset\Real$ define the probability measure
		\begin{equation}\label{eq:MC_EXTQ}
			\varPsi_n (B) = \left( \varGamma_n \otimes \varLambda_n \right) \left( \{ (\bfvartheta, \theta) \in \Theta_n\times(0,1] : \mu+\sigma Q_\gamma(\tau_E^{m\theta}) \in B \} \right).
		\end{equation} 
		Then, the $(1-\alpha)$-asymmetric confidence interval estimator for $Q_0(\tau_E)$, defined by the quantiles of $\varPsi_n$, satisfies
		\begin{equation}\label{eq:ACI_EXTQ}
			\Prob \left( Q_0(\tau_E) \in \left[ \varPsi_n ^{-1} (\alpha / 2); \varPsi_n ^{-1} (1 - \alpha / 2) \right]\right) = 1 - \alpha + o (1), \quad \forall\, \alpha \in (0, 1).
		\end{equation}
		\item Work under the conditions of Theorem \ref{thm:RETURNLV} with $\lambda_0 = 0$. For every measurable set $B\subset \Real$ define the probability measure
		\begin{equation}\label{eq:MC_RL}
			\varPsi_n (B) = \varGamma_n \left( \{ \bfvartheta \in \Theta_n : \mu+\sigma Q_\gamma(\tau^{m/m^{\star}}) \in B \} \right).
		\end{equation} 
		Then, the $(1-\alpha)$-asymmetric confidence interval estimator for $Q_0(\tau_E)$, defined by the quantiles of $\varPsi_n$, satisfies
		\begin{equation}\label{eq:ACI_RL}
			\Prob \left(R_0^{(m^\star)}(\tau) \in \left[ \varPsi_n^{-1} (\alpha / 2); \varPsi_n^{-1} (1 - \alpha / 2) \right]\right) = 1 - \alpha + o (1), \quad \forall\, \alpha \in (0, 1).
		\end{equation}
	\end{itemize}
\end{cor}
\begin{rem} \label{rem:asym_CI}
	Confidence intervals in \eqref{eq:ACI_EXTQ} and \eqref{eq:ACI_RL} can be approximated by a Monte Carlo approach. We generate replications $\bfvartheta^\star_i$ and $\theta^\star_i$ with $i=1,\ldots,M$, for a large $M\geq0$, from the three- and one-dimensional Gaussian distributions corresponding to the probability measure $\Gamma_n$ and $\Lambda_n$ in Corollary \ref{cor:CONFINT}. Then, for any $x\in\Real$, empirical estimates of the probability measures in \eqref{eq:MC_EXTQ} and \eqref{eq:MC_RL} are obtained respectively by  
\begin{align*}
\widehat{\varPsi}_n(x)&=\frac{1}{M} \sum_{i=1}^M\indic\left(\mu^\star_i+\sigma^\star_i Q_{\gamma^\star_i}\left(\tau_E^{m\theta^\star_i}\right)\leq x\right),\\
		\widehat{\varPsi}_n(x)&=\frac{1}{M}\sum_{i=1}^M\indic\left(\mu_i^\star+\sigma_i^\star Q_{\gamma^\star_i}\left(\tau^{m/m^{\star}}\right)\leq x\right).
\end{align*}
Then, for any $\alpha\in(0,1)$, the corresponding empirical quantiles $\widehat{\varPsi}_n^{-1}(\alpha/2)$ and $\widehat{\varPsi}_n^{-1}(1-\alpha/2)$ leads to the desired intervals.
As stated in Corollary \ref{cor:CONFINT},  for the theoretical result \eqref{eq:ACI_EXTQ} to hold, we only need an estimator $V_n$ of variance of $\widehat{\theta}_n$ that converges to a finite constant, see Section A.4 of the supplement for the motivation. For practical purposes, the selection of $V_n$  matters and in this regard we use the estimator in \eqref{eq:est_var_est_theta}.
\end{rem}
\begin{rem}\label{rem:quantbias}
Concluding when there is no asymptotic bias, 
	the results in Theorems \ref{thm:RETURNLV}--\ref{thm:EXTRQUANT} can be extended to time series models for which Condition \ref{cond:second_order}\ref{en_sec_order_1} is difficult to check (or not satisfied), by replacing such an assumption with Condition \eqref{eq:weaker} and further assuming that
	\begin{equation}\label{eq:quantbias}
		\frac{b( s \theta_0 m ) -b_{m,0}- a_{m,0} Q_{\gamma_0}(e^{-1/s}) }{a_{m,0} q_{\gamma_0}(e^{-1/s})}=o(1/\sqrt{k}),
	\end{equation}
	where either $s=(1+\epsilon_n)m^\star/(m (-\log \tau))$, for all $\epsilon_n$ such that $\epsilon_n =o(1/\sqrt{k})$, or  $s=1/(\theta_0 m (-\log \tau_E))$, for extreme return level or extreme VaR, respectively.
	See Section B for a formal discussion and Theorem 6.2 in  \cite{drees2003} for analogous results for the peaks over a threshold method.
\end{rem}

%
\section{Bayesian inference}\label{sec:Bayes_inf}
%

We propose a Bayesian procedure for inference on the extremal index $\theta_0$, GEV parameters $\bfvartheta_0$, return level $R^{(m^\star)}_0(\tau)$ and the extrapolation of the extreme VaR $Q_0(\tau_E)$ and establish its asymptotic properties.
Bayesian inference for the joint parameter vector $(\theta,\bfvartheta)$, whose individual loglikelihoods are  $\mathscr{L}_n(\bfvartheta)$ in \eqref{eq:logllik} and $\mathscr{L}_n(\theta)$ in  \eqref{eq:ext_ind_likelihood},  requires the derivation of their pseudo joint log-likelihood process. However, due to the dependence between the sequences $M_{i,m}$ for $i=1,\ldots,k$ and $\widehat{Y}_{i,\lm}$ for $i=1,\ldots,\lk$, the analytical form of the joint likelihood becomes intractable. To address this issue, we propose conducting inference on the individual parameters $\theta$ and $\bfvartheta$ separately. We justify this approach by showing that the posterior distributions of $\theta$ and $\bfvartheta$ possess key asymptotic properties that ensure their reliability. Moreover, although the ideal approach would involve deriving the posterior distribution of the extreme VaR from the joint posterior of $(\theta,\bfvartheta)$, we demonstrate that using the individual posterior distributions still yields a posterior for extreme VaR with desirable asymptotic behavior. In particular, the coverage probability of the resulting credible intervals converges to the nominal level (provided that there is no asymptotic bias).

We suggest the following two separate Bayesian methods. 
Concerning the estimation of $\bfvartheta$, its posterior distribution is defined according to an empirical Bayes approach, precisely using the pseudo loglikelihood  $\mathscr{L}_n(\bfvartheta)$ and by specifying a class of data-dependent prior distributions by their density function of the form
$
\pi(\bfvartheta)= \pi_{\text{sh}}^{(n)}(\gamma) \pi_{\text{loc}}^{(n)}
\left( \mu\right) \pi_{\text{sc}}^{(n)}
\left( \sigma\right)$,  $\bfvartheta \in \Theta$,
where for each $n=1,2,\ldots$, $\pi_{\text{sh}}^{(n)}$, $\pi_{\text{loc}}^{(n)}$ and $\pi_{\text{sc}}^{(n)}$ are prior densities on $\gamma$, $\mu$ and $\sigma$, whose expression may or may not depend on $n$. 
This formulation allows to elicit flexible and simple joint densities. We prescribe the minimal conditions under which the posterior satisfies favorable asymptotic properties that guarantee its accuracy. 
\begin{cond}\label{cond:prior_GEV_par}
	The densities $\pi_{\text{sh}}^{(n)}$, $ \pi_{\text{loc}}^{(n)}$ and $\pi_{\text{sc}}^{(n)}$ are such that with probability tending to one as $n \rightarrow \infty$:
	\begin{inparaenum}
		%
		\item\label{cond:pish} $\pi_{\text{sh}}^{(n)}$ is a positive function such that
		\begin{inparaenum}
			\item $\pi_{\text{sh}}^{(n)}$ is continuous and non-zero at $\gamma_0$,
			\item $\exists\, q < 1$ such that $\supp(\pi_{\text{sh}}^{(n)})\subset(-1/2, k^{q})$,
			\item $\int_{-1/2}^{k^{q}} \pi_{\text{sh}}^{(n)} (\gamma) \diff \gamma < \infty$.
		\end{inparaenum}
		\item\label{cond:piloc} $\pi_{\text{loc}}^{(n)}$ is a positive function such that
		\begin{inparaenum}
			\item $\int_{-\infty}^{\infty} \pi_{\text{loc}}^{(n)} (\mu) \diff \mu < \infty$,
			\item there is $\delta > 0$ such that $\pi_{\text{loc}}^{(n)} (b_{m,0}) a_{m,0} > \delta$ and for any $\eta > 0$ there is $\varepsilon > 0$ such that for all  $- \varepsilon \leq \tau \leq \varepsilon$
			\begin{equation*}
				(1 - \eta)  \pi_{\text{loc}}^{(n)} (b_{m,0})  \leq \pi_{\text{loc}}^{(n)} (a_{m,0}  \tau + b_{m,0}) \leq (1 + \eta) \pi_{\text{loc}}^{(n)} (b_{m,0}). 
			\end{equation*}
		\end{inparaenum}
		\item\label{cond:pisc} $\pi_{\text{sc}}^{(n)}$ is a positive function such that
		\begin{inparaenum}
			\item $\int_{0}^{\infty} \pi_{\text{sc}}^{(n)} (\sigma) \diff \sigma < \infty$,
			\item there is $\delta > 0$ such that $\pi_{\text{sc}}^{(n)} (a_{m,0}) a_{m,0} > \delta$ and for any $\eta > 0$ there is $\varepsilon > 0$ such that for all $1 - \varepsilon \leq \tau \leq 1 + \varepsilon$
			\begin{equation*}
				(1 - \eta)  \pi_{\text{sc}}^{(n)} (a_{m,0})  \leq \pi_{\text{sc}}^{(n)} (a_{m,0} \tau) \leq (1 + \eta) \pi_{\text{sc}}^{(n)} (a_{m,0}). 
			\end{equation*}%
		\end{inparaenum}
	\end{inparaenum}
\end{cond}
\begin{rem}\label{rem:prior}
A specific family of prior distributions satisfying Condition \ref{cond:prior_GEV_par} is the following class of data dependent priors. Consider some consistent estimators $\widehat{\sigma}_{n}$,  $\widehat{\mu}_{n}$ of $a_{m,0}, b_{m,0}$, in the sense that $\widehat{\sigma}_{n}/a_{m,0}=1+o_\Prob(1)$ and $(\widehat{\mu}_{n}-b_{m,0})/a_{m,0}=o_\Prob(1)$, e.g., the MLE in \eqref{eqMLE}.
Let $\pi_{\text{loc}}^{(n)} (\mu) = \pi_{\text{loc}} ((\mu - \widehat{\mu}_{n}) / \widehat{\sigma}_{n}) / \widehat{\sigma}_{n}$ and $\pi_{\text{sc}}^{(n)} (\sigma) = \pi_{\text{sc}}(\sigma / \widehat{\sigma}_{n}) / \widehat{\sigma}_{n}$.
with $\pi_{\text{loc}} : \mathbb{R} \rightarrow [0, \infty)$ and $\pi_{\text{sc}} : (0, \infty) \rightarrow [0, \infty)$ continuous functions such that $\mu \rightarrow \pi_{\text{loc}} (\mu u + v)$ is uniformly integrable for $u$ and $v$ close to $1$ and $0$, respectively, and $\sigma \rightarrow \pi_{\text{sc}} (\sigma u)$ is uniformly integrable for $u$ close to one.
Finally, let $\pi_{\text{sh}}$ be any continuous density function on $(-1/2, \infty)$ and let $\pi_{\text{sh}}^{(n)} (\gamma) \propto \pi_{\text{sh}} (\gamma) \ind(\gamma < k^{q})$ for some $q \in (0, 1)$.
\end{rem}

For all measurable sets $B\subset \Theta_n$, define the pseudo posterior distribution on $\bfvartheta$,
$$
\varPi_n(B)= \frac{\int_B \exp(k \mathscr{L}_n(\bfvartheta))\pi(\bfvartheta)\diff \bfvartheta}{\int_\Theta \exp(k \mathscr{L}_n(\bfvartheta))\pi(\bfvartheta)\diff \bfvartheta}.
$$
Note that $\varPi_n$ is not a genuine posterior distribution for several reasons. First, the pseudo log-likelihood $\mathscr{L}_n(\bfvartheta)$ is derived from a misspecified model. Second, the data $M_{i,m}$, $i=1,\ldots,k$, are not exactly independent, whereas $\mathscr{L}_n(\bfvartheta)$ corresponds to a likelihood for independent observations. Finally, and importantly, if a data-dependent prior is used, the data may effectively be used twice. Despite these limitations, the following result demonstrates that $\varPi_n$ still retains favorable asymptotic properties.
\begin{theo}[$\bfvartheta$-Posterior asymptotics]\label{theo:post_GEV_par}
	Under Conditions \ref{cond:second_order}, \ref{cond:ts_cond} and \ref{cond:prior_GEV_par} we have 
	\begin{itemize}
		\item (Contraction rate)  For any $C=C_n$ with $C\to\infty$ as $n\to\infty$ and $C=o(\sqrt{k})$ there is a $\varepsilon>0$ such that
		$$
		\varPi_n
		\left(
		\bfvartheta\in\Theta_n: (\gamma-\gamma_0)^2+\left(\frac{\mu-b_{m,0}}{a_{m,0}}\right)^2+\left(\frac{\sigma}{a_{m,0}}-1\right)^2 >\frac{C^2}{k}
		\right)
		=O_{\Prob}
		\left(
		e^{-\varepsilon C^2}
		\right);
		$$
		\item (Bernstein-von Mises) 
		For all measurable sets $B\subset \Theta_n$, as $n\to\infty$  
		$$
		\sup_{B\subset \Theta_n}
		\left|
		\varPi_n
		\left(
		\sqrt{k}\left(
		\gamma-\gamma_0, 
		\frac{\mu-b_{m,0}}{a_{m,0}},
		\frac{\sigma}{a_{m,0}}
		-1\right)\in B
		\right)-
		\mathcal{N}
		\left(B;
		\bfbeta ,\bfSigma
		\right)
		\right|
		=o_{\Prob}(1),
		$$
		where $\bfbeta= \sqrt{k}  \bfI^{-1}(\bftheta_0) \bfS_n(\bftheta_0)$ and $\bfSigma=\bfI^{-1}(\bftheta_0)$;
		\item (Coverage probability) For any $\alpha\in(0,1)$, let $\varPi^{-1}_{n,j}(\alpha)$, $j=1,2,3$, be the $\alpha$-quantile of the marginal posterior distribution $\varPi_{n,j}$ concerning the parameter $\gamma$, $\mu$ and $\sigma$, respectively. Let $I_{n,j}=[\varPi^{-1}_{n,j} (\alpha/2), \varPi^{-1}_{n,j}(1-\alpha/2)]$ be $(1-\alpha)$-asymmetric credible intervals. If $\lambda_0=0$,  then
		\begin{equation}\label{eq:credible_int_GEV}
			\lim_{n \to \infty} \Prob( \gamma_0 \in I_{n,1} )=
			\lim_{n \to \infty} \Prob( b_{m,0} \in I_{n,2} )=
			\lim_{n \to \infty} \Prob( a_{m,0}  \in I_{n,3} )=
			1-\alpha.
		\end{equation}
	\end{itemize}
\end{theo}
Leveraging on the Bernstein-von Mises result of Theorem \eqref{theo:post_GEV_par}, an approximation of an asymptotic $(1-\alpha)$-equi-tailed credible interval estimator for $\gamma_0$  is given by 
\begin{equation}\label{eq:BCI_gamma}
	\left[
	\overline{\gamma}_n-\widehat{b}_n\pm z_{1-\alpha/2}\overline{\psi}_n
	\right],
\end{equation}
where $\overline{\gamma}_n$ is the posterior mean, $\widehat{b}_n$ is an estimator of $b_0$ the first element of $\sqrt{k}\lambda_0 \theta_{0}^{\rho_0} \bfI^{-1}(\bftheta_0) \bfb_0$ and $\overline{\psi}_n$ is the posterior standard deviation, here seen as an estimator of $\psi_0$, i.e. the square root of the first diagonal element of $k^{-1} \bfI^{-1}(\bftheta_0)$.

Now, let $\bfvartheta$ follow the distribution $\varPi_n$. Then, for every suitable measurable sets $B\subset\Real$ and $\tau\in(0,1)$, we define the posterior distribution of the $T$-return level corresponding to the sample maximum $M_{m^\star}$ of block size $m^\star$, induced by the transformation $\widetilde{R}_{\bfvartheta}(\tau)=\mu +\sigma Q_{\gamma}(\tau^{m^\star / m})$, as
\begin{equation}\label{eq:posterior_ret_lev}
	\varPi_n^{(RL)}(B)= \Pi_n \left( \left\{ \bfvartheta \in \Theta_n : \widetilde{R}_{\bfvartheta}(\tau) \in B \right\}\right).
\end{equation}
\begin{cor}[$T$-Return level posterior asymptotics]\label{cor:POSTERIORRETURNLV}
	Under the conditions in Theroem \ref{thm:RETURNLV} and Theorem \ref{theo:post_GEV_par} we have 
	\begin{itemize}
		\item{(Contraction rate)} For any sequence $C_n$ such that $C_n \rightarrow \infty$ as $n \rightarrow \infty$ and $ -C_n\log( -m/m^\star\log\tau)= o(\sqrt{k})$,  there is a $\varepsilon > 0$ such that 
		\begin{equation*}
			\varPi_n^{(RL)}\left( \widetilde{R}_{\bfvartheta}(\tau) \in \Real : \biggr\vert \frac{\widetilde{R}_{\bfvartheta}(\tau) - R_0^{(m^*)}(\tau)}{a_{m,0} q_{\gamma_0} (\tau^{m/m^\star})} \biggr\vert > \frac{C_n}{\sqrt{k}} \right) = O_{\Prob} \left( e^{-\varepsilon C_{n}^2}\right)
		\end{equation*}
		\item{(Bernstein-von Mises)} For all measurable sets $B \in \mathbb{R}$ we have, as $n\to\infty$
		\begin{equation*}
			\underset{B \in \Real}{\sup} \biggr\vert \varPi_n^{(RL)} \left(  \widetilde{R}_{\bfvartheta}(\tau) \in \Real : 
			\frac{\sqrt{k}(\widetilde{R}_{\bfvartheta}(\tau) - R_0^{(m^*)}(\tau))}{a_{m,0} q_{\gamma_0} (\tau^{m/m^\star})}\in B\right) 
			-\mathcal{N} \left(B; \beta, \Sigma \right) \biggr\vert = o_{\Prob} (1),
		\end{equation*}
where 
$\beta \overset{d}{\rightarrow} \mathcal{N} (\lambda_0\theta_{0}^{\rho_{0}} \mathbf{v}_{\gamma_0}^{\top} \bfI^{-1}(\bftheta_0) \bfb_0, \bfSigma)$ and $\Sigma = \mathbf{v}_{\gamma_0}^{\top} \bfI^{-1}(\bftheta_0)  \mathbf{v}_{\gamma_0}$ with $\bfb_0$ defined in Section D of the Supplement and $\mathbf{v}_{\gamma_0}$ in Theorem \ref{thm:RETURNLV}.
		\item{(Coverage probability)} For any $\alpha \in (0, 1)$, let $ \varPi_n^{(RL)-1} (\alpha)$ be the $\alpha$-quantile of $ \varPi_n^{(RL)}$. If $\lambda_0 = 0$, then
\begin{equation}\label{eq:ABCI_RL}
\underset{n \rightarrow \infty}{\lim} \Prob \left(R_0^{(m^*)}(\tau) \in [\varPi_n^{(RL)-1}(\alpha / 2); \varPi_n^{(RL)-1} (1 - \alpha / 2)]  \right)=1 - \alpha.
\end{equation} 
\end{itemize}
\end{cor}
On the basis of the  BvM result of Theorem \ref{eq:posterior_ret_lev}, an approximation of an asymptotic $(1-\alpha)$-equi-tailed credible interval for $R^{(m^\star)}_{0}(\tau)$ is given by
\begin{equation}\label{eq:SBCI_RL}
	\left[
	\overline{R}_n(\tau^{m / m^\star}) -  \widehat{b}_n \pm z_{1-\alpha/2}\overline{\psi}_n
	\right],
\end{equation}
where $\overline{R}_n(\tau^{m / m^\star})$ is the return level posterior mean, $\widehat{b}_n$ is an estimator of the true bias  $b_0=a_{m,0}q_{\gamma_0}(\tau^{m/m^\star}) \sqrt{k}\lambda_0\theta_{0}^{\rho_{0}} \left(\mathbf{v}_{\gamma_0}^{\top} \bfI^{-1}(\bftheta_0) \bfb_0 + b_{\gamma_0}\right)$, $\overline{\psi}_n$ is the return level posterior standard deviation, here seen as an estimator of $\psi_0=\sqrt{k}  a_{m,0} q_{\gamma_0} (\tau^{m / m^{\star}} ) \sqrt{ \Sigma }$.

We next focus now on the Bayesian procedure for the estimation of $\theta_0\in(0,1]$. Its posterior distribution is defined using the loglikelihood  $\mathscr{L}_n(\theta)$ and by specifying a prior distribution as a mixture for all measurable sets $B \subset [0,1]$
\begin{equation}\label{eq:prior_theta}
	\Pi(B)= (1-p)\widetilde{\Pi}(B \cap (0,1) )+ p\delta_{1}(B),
\end{equation}
where $\widetilde{\Pi}$ is an informative proper prior distribution with continuous density $\widetilde{\pi}$ with support $(0,1)$, e.g. a Beta, uniform, truncated, etc., $\delta_1$ is the Dirac measure and $p\in [0,1)$ is a point mass on $1$. For all measurable sets $B\subset (0,1]$, we define the following pseudo posterior distribution for $\theta$,
$$
\varPhi_n(B)= \frac{\int_B \exp(\lk \mathscr{L}_n(\theta))\Pi(\diff\theta)}{\int_{(0,1]} \exp(\lk \mathscr{L}_n(\theta))\Pi (\diff \theta)}.
$$
The posterior distribution $\varPhi_n$ is not a genuine posterior in the strict Bayesian sense, because the exponential log-likelihood $\mathscr{L}_n(\theta)$ , defined in \eqref{eq:ext_ind_likelihood}, is constructed from the variables  $\widehat{Y}_{i,\lm}$, $i=1,\ldots,\lk$,  which serve as proxies for the true $Y_{i,\lm}$ ones, which  are only asymptotically exponentially distributed. Then, the exponential model is misspecified. Nevertheless, we can still show the following good asymptotic properties for $\varPhi_n$.
\begin{theo}[$\theta$-Posterior asymptotics]\label{theo:post_extremal_index}
	Under Condition 2.1 in \cite{berghaus2018} we have 
	\begin{itemize}
		\item (Contraction rate)  For any $C=C_n$ with $C\to\infty$ as $n\to\infty$ and $C=o(\lk^{1/2})$ there is a $\varepsilon>0$ such that
		$$
		\varPhi_n
		\left(
		\theta\in(0,1]: |\theta-\theta_0| >\frac{C}{\lk^{1/2}}
		\right)
		=O_{\Prob}
		\left(
		e^{-\varepsilon C^2}
		\right);
		$$
		Moreover, if $\theta_0 = 1$ then $\varPhi_{n} \left( \theta \neq \theta_0 \right) = o_{\Prob} (1)$.
		\item (Bernstein-von Mises) 
		If $\theta_0\in(0,1)$, for all measurable sets $B\subset \Real$ we have, as $n\to\infty$ 
		$$
		\sup_{B\subset (0,1)}
		\left|
		\varPhi_n
		\left(
		\lk^{1/2}\left(
		\theta-\theta_0\right)\in B
		\right)-
		\mathcal{N}
		\left(B;
		\beta,\theta_{0}^{2}
		\right)
		\right|
		=o_{\Prob}(1),
		$$
		where $\beta \overset{d}{\rightarrow} \mathcal{N} (0, \theta_{0}^4 \sigma^2(\theta_0))$ and $0 < \sigma^{2}(\theta_0) \neq \theta_{0}^{-2}$ is explicitly given in Theorem 3.2 in \cite{berghaus2018}.
	\end{itemize}
\end{theo}
Note that in Theorem \ref{theo:post_extremal_index} the asymptotic variance of the posterior distribution $\varPhi_n$ is equal to $\theta_0^2$ and does not coincide with the asymptotic variance of posterior random center $\beta$, which is instead equal to $\theta_{0}^4 \sigma^2(\theta_0)$.  The effect of this is that the coverage probability of credible intervals for $\theta_0$ derived from $\varPhi_n$ can not asymptotically achieve the nominal level, in contrast to their frequentist counterparts based on the MLE. To overcome this limitation, we propose an adjustment to the posterior distribution given by a reparametrization of the log-likelihood by the linear transformation
$$
\theta \mapsto \widehat{\theta}_n + (\theta-\widehat{\theta}_n)(\widehat{\theta}_n\widetilde{\sigma}_{n})^{-1},
$$
where $ \widehat{\theta}_n$ is the MLE of $\theta_0$ in \eqref{eq:ext_ind_MLE} and  $\widetilde{\sigma}_{n}^{2}$ is an estimator of $\sigma^2(\theta_0)$ satisfying $\widetilde{\sigma}_{n}^{2}=\sigma^2(\theta_0) + o_\Prob(1)$.
This correction is designed to suitably rescale the variability of the resulting posterior distribution. We refer to this as the quadrature adjustment, and based on it, we define the corresponding quadrature-adjusted posterior distribution
\begin{equation}\label{eq:adjusted_posterior}
	\varPhi_{n}^{(A)} (B) = \frac{\int_{B} \exp \left( \lk \mathscr{L}_{n} \left(\widehat{\theta}_n + (\theta-\widehat{\theta}_n)(\widehat{\theta}_n\widetilde{\sigma}_{n})^{-1}\right) \right) \Pi (\diff \theta)}{\int_{(0, 1]} \exp \left( \lk \mathscr{L}_{n} \left(\widehat{\theta}_n + (\theta-\widehat{\theta}_n)(\widehat{\theta}_n\widetilde{\sigma}_{n})^{-1}\right) \right) \Pi (\diff \theta)}.
\end{equation}

\begin{theo}[Adjusted $\theta$-Posterior asymptotics]\label{theo:post_extremal_index_adj}
	Under Condition 2.1 in \cite{berghaus2018} we have 
	\begin{itemize}
		\item (Contraction rate)  For any $C=C_n$ with $C\to\infty$ as $n\to\infty$ and $C=o(\lk^{-1/2})$ there is a $\varepsilon>0$ such that
		$$
		\varPhi_{n}^{(A)}
		\left(
		\theta\in(0,1]: |\theta-\theta_0| >\frac{C}{\lk^{1/2}}
		\right)
		=O_{\Prob}
		\left(
		e^{-\varepsilon C^2}
		\right);
		$$
		Moreover, if $\theta_0 = 1$ then $\varPhi_{n}^{(A)} \left( \theta \neq \theta_0 \right) = o_{\Prob} (1)$.
		\item (Bernstein-von Mises) 
		If $\theta_0\in(0,1)$, for all measurable sets $B\subset \Real$ we have, as $n\to\infty$ 
		$$
		\sup_{B\subset (0,1)}
		\left|
		\varPhi_{n}^{(A)}
		\left(
		\lk^{1/2}\left(
		\theta-\theta_0\right)\in B
		\right)-
		\mathcal{N}
		\left(B;
		\beta,\theta_{0}^{4} \sigma^2(\theta_0)
		\right)
		\right|
		=o_{\Prob}(1),
		$$
		where $\beta \overset{d}{\rightarrow} \mathcal{N} (0, \theta_{0}^4 \sigma^2(\theta_0))$ and $\sigma^{2}(\theta_0)$ is as in Theorem \ref{theo:post_extremal_index};
		\item (Coverage probability) For any $\alpha\in(0,1)$, let $\varPhi^{(A)-1}(\alpha)$ be the $\alpha$-quantile of $\varPhi^{(A)}$. If  $\theta_0\in(0,1)$ and $\lambda_0=0$, then
		\begin{equation}\label{eq:ABCI_EXT_IND}
			\lim_{n \to \infty} \Prob\left( \theta_0 \in \left[ \varPhi^{(A)-1}(\alpha/2); \varPhi^{(A)-1} (1-\alpha/2) \right]\right)=
			1-\alpha.
		\end{equation}
	\end{itemize}
\end{theo}
Exploiting the  BvM result of Theorem \ref{theo:post_extremal_index_adj}, an approximation of an asymptotic $(1-\alpha)$-equi-tailed credible interval for $\theta_{0}(\tau)$ is given by
\begin{equation}\label{eq:SBCI_EXT_IND}
	\left[
	\overline{\theta}_n\pm z_{1-\alpha/2}\overline{\psi}_n
	\right],
\end{equation}
where $\overline{\theta}_n$ is the extremal index posterior mean and $\overline{\psi}_n$ is the extremal index posterior standard deviation.

Now, let $\bfvartheta$ and $\theta$ follow the distribution $\varPi_n$ and $\varPhi_n$, respectively. Then, for every suitable measurable sets $B\subset\Real$ and extreme level $\tau_E\to 1$ as $n\to\infty$, we define the posterior distribution of the extreme VaR, induced by the transformation $\widetilde{Q}_{\bfvartheta,\theta}(\tau_E)=\mu +\sigma Q_{\gamma}(\tau_E^{m\theta})$, as
\begin{equation}\label{eq:posterior_ext_VaR}
	\varPhi_n^{(VaR)}(B)= (\varPi_n \otimes\varPhi_n^{(A)}) \left( \left\{ (\bfvartheta, \theta) \in \Theta_n \times (0, 1] : \widetilde{Q}_{\bfvartheta,\theta}(\tau_E) \in B \right\}\right).
\end{equation}
\begin{cor} \label{cor:POSTERIOREXTRQUANT}
	Suppose that $\lim_{n \rightarrow \infty} k / \lk = c \in [0, \infty)$.
	Then, under the conditions in Theroem \ref{thm:EXTRQUANT}, Theorem \ref{theo:post_GEV_par} and Theorem \ref{theo:post_extremal_index_adj} we have
	\begin{itemize}
		\item{(Contraction rate)} For any $C=C_n$ with $C \rightarrow \infty$ as $n \rightarrow \infty$ and $-C\log( -\theta_0 m\log\tau_E) = o (\sqrt{k})$, there is a $\varepsilon > 0$ such that 
		\begin{equation*}
			\varPhi_n^{(VaR)} \left( \widetilde{Q}_{\bfvartheta,\theta}(\tau_E) \in \Real : \biggr\vert \frac{\widetilde{Q}_{\bfvartheta,\theta}(\tau_E) - Q_0(\tau_E)}{a_{m,0} q_{\gamma_0} (\tau_E^{m\theta_0})} \biggr\vert > \frac{C}{\sqrt{k}} \right) = O_{\Prob} \left( e^{-\varepsilon C^2}\right)
		\end{equation*}
		\item{(Bernstein-von Mises)} For all measurable sets $B \in \mathbb{R}$ we have, as $n\to\infty$
		\begin{equation*}
			\underset{B \in \Real}{\sup} \biggr\vert \varPhi_n^{(VaR)}\left( \widetilde{Q}_{\bfvartheta,\theta}(\tau_E) \in \Real : \frac{\sqrt{k}(\widetilde{Q}_{\bfvartheta,\theta}(\tau_E) - Q_0(\tau_E))}{a_{m,0} q_{\gamma_0} (\tau_E^{m\theta_0})} \in B  \right) - 
			\mathcal{N} \left(B; \beta, \Sigma \right) \biggr\vert = o_{\Prob} (1),
		\end{equation*}
		where 
		$\beta \overset{d}{\rightarrow} \mathcal{N} (\lambda_0\theta_{0}^{\rho_{0}}\left(\mathbf{v}_{\gamma_0}^{\top} \bfI^{-1}(\bftheta_0) \bfb_0+b_{\gamma_0}\right), \bfSigma)$, $\Sigma = \mathbf{v}_{\gamma_0}^{\top} \bfI^{-1}(\bftheta_0)  \mathbf{v}_{\gamma_0}$ with $\bfb_0$ defined in Section D, $\mathbf{v}_{\gamma_0}$ and $b_{\gamma_0}$ are defined as in Theorem \ref{thm:RETURNLV} and \ref{thm:EXTRQUANT}, respectively.
		\item{(Coverage probability)} For any $\alpha \in (0, 1)$, let $\varPhi_n^{(VaR)-1} (\alpha)$ be the $\alpha$-quantile of $\varPhi_n^{(VaR)}$. If $\lambda_0 = 0$, then
		\begin{equation}\label{eq:ABCI_EXT_VAR}
			\underset{n \rightarrow \infty}{\lim} \ \Prob \left(Q_0(\tau_E) \in  [\varPhi_n^{(VaR)-1}(\alpha / 2); 
			\varPhi_n^{(VaR)-1}(1 - \alpha / 2)]\right)= 1 - \alpha.
		\end{equation} 
	\end{itemize}
\end{cor}
\begin{rem}\label{rem:details_post_asym}
	The results in Corollary \ref{cor:POSTERIOREXTRQUANT} still hold if we use the unadjusted posterior $\varPhi_n$ instead of the adjusted version $\varPhi^{(A)}_n$, in formula \eqref{eq:posterior_ext_VaR}. 
\end{rem}
On the basis of the BvM result of Corollary \ref{cor:POSTERIOREXTRQUANT}, an approximation of an asymptotic $(1-\alpha)$-equi-tailed credible interval for $Q_{0}(\tau_E)$ is given by
\begin{equation}\label{eq:SBCI_EXT_VAR}
	\left[
	\overline{Q}_n(\tau_E) -  \widehat{b}_n \pm z_{1-\alpha/2}\overline{\psi}_n
	\right],
\end{equation}
where $\widehat{b}_n$ is an estimator of $b_0=a_{m,0}q_{\gamma_0}(\tau_E^{m \theta_0})\sqrt{k} \lambda_0\theta_{0}^{\rho_{0}}\left(\mathbf{v}_{\gamma_0}^{\top} \bfI^{-1}(\bftheta_0) \bfb_0+b_{\gamma_0}\right)$, $\overline{Q}_n(\tau_E)$ is the extreme VaR posterior mean and $\overline{\psi}_n$ is the extreme VaR posterior standard deviation, here seen an estimator of $\psi_0=\sqrt{k}  a_{m,0} q_{\gamma_0} (\tau^{m / m^{\star}} ) \sqrt{ \Sigma }$.

\begin{rem}
Analogous considerations to those discussed in Remark \ref{rem:quantbias} for the frequentist framework also apply in the Bayesian setting. In particular, the posterior contraction rates and asymptotic normality for $\boldsymbol{\vartheta}$ established in Theorem \ref{theo:post_GEV_par} remain valid when assumption \eqref{eq:weaker} is used in place of Condition \ref{cond:second_order}\ref{en_sec_order_1}. Moreover, the results on posterior inference for return levels and extreme VaR, as stated in Corollaries \ref{cor:POSTERIORRETURNLV}--\ref{cor:POSTERIOREXTRQUANT}, continue to hold under this weaker condition, provided that assumption \eqref{eq:quantbias} is additionally imposed in the extrapolation regime. Further details are provided in Section B of the Supplement.
\end{rem}

%
%
\section{Numerical experiments}\label{sec:simulations}
%
%
We present a simulation study to evaluate the finite-sample performance of the ML estimators introduced in Section \ref{sec:LAT} and the Bayesian procedures developed in Section \ref{sec:Bayes_inf}. Additional simulation details and results are provided in Section E of the supplement.
Our aim is assessing the frequentist coverage of asymptotic $(1 - \alpha)$ confidence intervals—both equi-tailed (symmetric) and asymmetric—and their Bayesian counterparts for key quantities discussed in the paper.
We also compare the point estimators derived from posterior medians with the ML estimators.

Simulations are made using models from \cite{berghaus2018}, including:
\begin{enumerate}
	\item[a)] the ARMAX model $X_{t+1} = \max(\eta X_{t}, (1 - \eta) Z_{t+1})$, where the innovations $Z_t$ are independent with unit-Fr\'{e}chet distribution $\Prob(X_t\leq x)=e^{1/x}$ and $\eta=(0, 0.5)$;
	\item[b)] the Markovian Copula models $X_{t+1} = F^{-1} (U_{t+1})$, with $(1 - U_{t}, 1 - U_{t+1})$ following the Clayton Copula  $C_{\eta} (u_1, u_2) = (u_{1}^{-\eta} + u_{2}^{-\eta} - 1)^{-1 / \eta}$, where $\eta =(0.41, 1.06)$, and $F(x)=1-e^{-x}$ for $x\geq0$ in one case and $F(x)=1-(1-x)^3/9$ in the other.
	\item[c)] the ARCH model $X_{t+1} =(2\times10^{-5} + \eta X_{t}^2)^{-1/2} Z_{t+1}$, where $\eta=(0.5, 0.99)$ and $Z_t$ are independent standard normal innovations.
\end{enumerate}
In the first model, the one-dimensional stationary distribution is a member of the GEV class, with $\gamma_0=1$. In this case, $a_{m,0}=b_{m,0}= \theta_0 m$ and there is no need for the second-order assumptions in Condition \ref{cond:second_order} for the theory to work. 
The values $\eta = (0, 0.5)$ correspond to the cases of perfect independence and strong dependence, respectively, and accordingly we have  $\theta_0 = (1, 0.5)$, respectively, see Section 6.1 in \cite{berghaus2018} for details.
In the second and third models, the one-dimensional stationary distribution are unit-exponential and a specific member of the power-law class $F(x)=1-K(x^\star-x)^\nu$ with the end-point $x^\star$, shape parameter $\nu>0$ and constant $K>0$ \citep[][Example 3.3.16]{embrechts1997modelling}, and therefore they are in the domain of attraction of a Gumbel family or GEV with $\gamma_0=0$ and Reverse-Weilbull family or GEV with $\gamma_0 = - 1 / 3$, respectively. For both distributions, \cite{padoan2024} showed that Condition \ref{cond:second_order}\ref{en_sec_order_1} is satisfied for $a_{m,0}= \theta_0 m  b_{0}^{'} (\theta_0 m)$, where $b_{0}^{'}(x)=(\partial / \partial x) b_0(x)$.
The values $\eta = (0.41, 1.06)$ correspond to the cases of weak dependence and strong dependence, respectively, and accordingly we have $\theta_0 \approx (0.80, 0.40)$, see \cite{perfekt1994} for the details on the computation of the extremal index. In the fourth model, the marginal distribution is heavy-tailed, and for values of $\alpha=(0.5,0.99)$, we obtain $\gamma_0\approx(0.211,0.493)$, as derived from \citet[][Theorem 2.1]{mikosch2000limit}. Moreover, the corresponding values of the extremal index are approximately $\theta_0\approx(0.832,0.565)$, representing the case of weak and strong dependence. In this setting, our theoretical results hold under the weaker condition \eqref{eq:weaker} discussed in Remark \ref{rem:weak_SOC}, which relaxes the standard second-order condition, and the bias assumption \eqref{eq:quantbias}.

We simulate time series with sample sizes $n \in {360, 1800, 3600, 10800}$ (representing $1, 5, 10$, and $30$ years of daily data). Disjoint BM are computed using block sizes $m \in {7, 30, 90, 180, 360}$. For inference, GEV parameters are estimated via MLE and Bayesian methods. Without loss of generality, we define the parameter space as  $\Theta_n=(-1/2, \sqrt{k}) \times \Real \times (0, \infty)$. For simplicity, we always define the small block $l=0$. We estimate $\bfvartheta_0$ with the MLE $\widehat{\bfvartheta}_n$ in \eqref{eqMLE} and 
$\theta_0$ with the MLE $\widehat{\theta}_n$ in \eqref{eq:ext_ind_MLE} and  its variance by the estimator $\widehat{\theta}_n^4\widetilde{\sigma}_{n}^{2}$, with $\widetilde{\sigma}_{n}^{2}$ is given in \eqref{eq:est_var_est_theta} and where we set $K=10$. The Bayesian posterior is computed using a hybrid of self-normalized importance sampling \citep{tokdar2010} and an adaptive random-walk Metropolis-Hastings algorithm \citep[][Section 4]{padoan2024}. A truncated Student-$t$ prior on the interval $(-1/2,  \sqrt{k})$ is used for $\gamma$, and the prior for $\theta$ is given by \eqref{eq:prior_theta}, combining a uniform density on $(0,1)$ with a point mass at 1 ($p = 0.1$). We set for simplicity $\tilde{m}=m$ and $\tilde{k}=k$. The adjusted posterior for $\theta$, obtained by formula \eqref{eq:adjusted_posterior}, is known in closed form apart from the normalizing constant which is evaluated via numerical integration, see Section E of the supplement for all the details. All the Bayesian procedures are available in the \textsf{R}-package \textbf{ExtremeRisks} \cite{ExtremeRisks}.

\begin{table}[t!]
	{\scriptsize
		\centering
		\caption{Coverage (strong dependence)}
		\resizebox{\textwidth}{!}{\begin{tabular}{c c c c|c c c c | c c c c | c c c c}
				\hline
				\multicolumn{4}{c}{} & \multicolumn{4}{c}{} & \multicolumn{8}{c}{Markov copula model ($\eta = 1.06$)} \\
				\multicolumn{4}{c}{} & \multicolumn{4}{c}{ARMAX ($\eta = 0.5$)} & \multicolumn{4}{c}{Exponential distr.} &  \multicolumn{4}{c}{Powerlaw distr.} \\
				\hline
				$n$ & $k$ & $m$ & Type 
				& $\gamma_0$ & $\theta_0$ &  RL & EQ 
				& $\gamma_0$ & $\theta_0$ &  RL & EQ
				& $\gamma_0$ & $\theta_0$ &  RL & EQ \\
				\hline
				360 & 12 & 30 & BS & 0.90 & 0.85 & 0.93 & 0.88 & 0.95 & 0.74 & 0.94 & 0.90 & 0.95 & 0.74 & 0.94 & 0.91 \\
				&  & & BA & 0.92 & 0.84 & 0.91 & 0.87 & 0.91 & 0.68 & 0.91 & 0.88 & 0.91 & 0.68 & 0.89 & 0.89 \\
				&  & & FS & 0.91 & 0.92 & 0.79 & 0.77 & 0.96 & 0.86 & 0.76 & 0.71 & 0.98 & 0.86 & 0.87 & 0.77 \\
				&  & & FA & - & - & 0.92 & 0.88 & - & - & 0.89 & 0.81 & - & - & 0.96 & 0.87 \\
				\hline
				& 51 & 7 & BS & 0.92 & 0.76 & 0.95 & 0.90 & 0.87 & 0.06 & 0.97 & 0.94 & 0.74 & 0.05 & 0.96 & 0.93 \\
				&  & & BA & 0.92 & 0.71 & 0.91 & 0.89 & 0.82 & 0.04 & 0.88 & 0.90 & 0.69 & 0.04 & 0.87 & 0.89 \\
				&  & & FS & 0.92 & 0.81 & 0.85 & 0.82 & 0.89 & 0.08 & 0.92 & 0.87 & 0.84 & 0.08 & 0.93 & 0.86 \\
				&  & & FA & - & - & 0.92 & 0.90 & - & - & 0.94 & 0.90 & - & - & 0.95 & 0.90 \\
				\hline
				1800 & 20 & 90 & BS & 0.92 & 0.91 & 0.94 & 0.90 & 0.94 & 0.87 & 0.93 & 0.90 & 0.95 & 0.87 & 0.94 & 0.92 \\
				&  & & BA & 0.93 & 0.89 & 0.92 & 0.89 & 0.92 & 0.82 & 0.91 & 0.89 & 0.94 & 0.83 & 0.90 & 0.89 \\
				&  & & FS & 0.93 & 0.93 & 0.81 & 0.79 & 0.95 & 0.91 & 0.81 & 0.78 & 0.99 & 0.91 & 0.86 & 0.81 \\
				&  & & FA & - & - & 0.94 & 0.91 & - & - & 0.90 & 0.85 & - & - & 0.95 & 0.88 \\
				\hline
				& 60 & 30 & BS & 0.94 & 0.92 & 0.95 & 0.92 & 0.94 & 0.66 & 0.96 & 0.94 & 0.91 & 0.66 & 0.94 & 0.92 \\
				&  & & BA & 0.94 & 0.91 & 0.94 & 0.91 & 0.91 & 0.60 & 0.93 & 0.93 & 0.89 & 0.60 & 0.91 & 0.90 \\
				&  & & FS & 0.94 & 0.94 & 0.86 & 0.84 & 0.95 & 0.71 & 0.90 & 0.87 & 0.96 & 0.71 & 0.90 & 0.86 \\
				&  & & FA & - & - & 0.95 & 0.92 & - & - & 0.94 & 0.91 & - & - & 0.94 & 0.90 \\
				\hline
				3600 & 20 & 180 & BS & 0.92 & 0.90 & 0.94 & 0.90 & 0.93 & 0.89 & 0.93 & 0.90 & 0.94 & 0.89 & 0.94 & 0.93 \\
				&  & & BA & 0.94 & 0.89 & 0.94 & 0.92 & 0.91 & 0.85 & 0.91 & 0.89 & 0.94 & 0.85 & 0.91 & 0.90 \\
				&  & & FS & 0.95 & 0.93 & 0.85 & 0.84 & 0.93 & 0.92 & 0.80 & 0.77 & 0.98 & 0.92 & 0.85 & 0.81 \\
				&  & & FA & - & - & 0.95 & 0.93 & - & - & 0.89 & 0.85 & - & - & 0.95 & 0.88 \\
				\hline
				& 40 & 90 & BS & 0.94 & 0.93 & 0.95 & 0.92 & 0.95 & 0.88 & 0.94 & 0.92 & 0.93 & 0.88 & 0.92 & 0.90 \\
				&  & & BA & 0.94 & 0.92 & 0.94 & 0.92 & 0.93 & 0.84 & 0.93 & 0.92 & 0.92 & 0.84 & 0.89 & 0.88 \\
				&  & & FS & 0.95 & 0.94 & 0.85 & 0.84 & 0.95 & 0.90 & 0.85 & 0.83 & 0.99 & 0.90 & 0.85 & 0.82 \\
				&  & & FA & - & - & 0.95 & 0.93 & - & - & 0.92 & 0.88 & - & - & 0.92 & 0.88 \\
				\hline
				& 120 & 30 & BS & 0.94 & 0.92 & 0.94 & 0.92 & 0.92 & 0.47 & 0.98 & 0.96 & 0.85 & 0.47 & 0.96 & 0.94 \\
				&  & & BA & 0.94 & 0.91 & 0.94 & 0.92 & 0.89 & 0.44 & 0.92 & 0.93 & 0.81 & 0.43 & 0.91 & 0.92 \\
				&  & & FS & 0.94 & 0.93 & 0.88 & 0.87 & 0.93 & 0.51 & 0.94 & 0.91 & 0.91 & 0.51 & 0.93 & 0.90 \\
				&  & & FA & - & - & 0.94 & 0.92 & - & - & 0.95 & 0.93 & - & - & 0.95 & 0.92 \\
				\hline
				10800 & 30 & 360 & BS & 0.93 & 0.93 & 0.95 & 0.91 & 0.94 & 0.92 & 0.93 & 0.91 & 0.94 & 0.92 & 0.94 & 0.93 \\
				&  & & BA & 0.94 & 0.91 & 0.94 & 0.91 & 0.92 & 0.89 & 0.92 & 0.90 & 0.93 & 0.89 & 0.91 & 0.90 \\
				&  & & FS & 0.94 & 0.94 & 0.84 & 0.83 & 0.93 & 0.93 & 0.83 & 0.81 & 0.95 & 0.93 & 0.85 & 0.83 \\
				&  & & FA & - & - & 0.95 & 0.93 & - & - & 0.90 & 0.87 & - & - & 0.93 & 0.89 \\
				\hline
				& 60 & 180 & BS & 0.94 & 0.94 & 0.95 & 0.93 & 0.95 & 0.90 & 0.95 & 0.94 & 0.93 & 0.90 & 0.93 & 0.91 \\
				&  & & BA & 0.94 & 0.93 & 0.94 & 0.92 & 0.94 & 0.88 & 0.94 & 0.93 & 0.93 & 0.88 & 0.91 & 0.89 \\
				&  & & FS & 0.94 & 0.95 & 0.86 & 0.86 & 0.95 & 0.92 & 0.88 & 0.86 & 0.96 & 0.92 & 0.85 & 0.84 \\
				&  & & FA & - & - & 0.95 & 0.93 & - & - & 0.93 & 0.90 & - & - & 0.90 & 0.88 \\
				\hline
				& 120 & 90 & BS & 0.95 & 0.94 & 0.95 & 0.93 & 0.95 & 0.81 & 0.96 & 0.95 & 0.92 & 0.82 & 0.94 & 0.93 \\
				&  & & BA & 0.95 & 0.93 & 0.95 & 0.93 & 0.94 & 0.78 & 0.94 & 0.93 & 0.91 & 0.78 & 0.92 & 0.92 \\
				&  & & FS & 0.95 & 0.94 & 0.88 & 0.88 & 0.95 & 0.84 & 0.92 & 0.90 & 0.96 & 0.84 & 0.90 & 0.88 \\
				&  & & FA & - & - & 0.95 & 0.93 & - & - & 0.94 & 0.92 & - & - & 0.93 & 0.91 \\
				\hline
				& 360 & 30 & BS & 0.95 & 0.88 & 0.95 & 0.93 & 0.81 & 0.08 & 0.98 & 0.96 & 0.62 & 0.08 & 0.97 & 0.97 \\
				&  & & BA & 0.95 & 0.87 & 0.95 & 0.93 & 0.89 & 0.07 & 0.89 & 0.92 & 0.57 & 0.07 & 0.89 & 0.93 \\
				&  & & FS & 0.95 & 0.90 & 0.92 & 0.90 & 0.84 & 0.09 & 0.97 & 0.95 & 0.70 & 0.09 & 0.97 & 0.96 \\
				&  & & FA & - & - & 0.94 & 0.93 & - & - & 0.93 & 0.94 & - & - & 0.96 & 0.96 \\
				\hline
		\end{tabular}}
		\label{table:covstrongdep} 
	}
	\footnotesize{
		Empirical coverage of symmetric (BS) and asymmetric (AS) 95\% Bayesian credible intervals and their frequentist counterparts (FS, FA), for $\gamma_0$, $\theta_0$, $R_{0}^{(n)}(0.9)$ (RL), and $Q_{0}(1 - 1 / n)$ (EQ), obtained with the ARMAX and Markov copula models and different settings: $n$, $k$ and $m$.}
\end{table}

Table \ref{table:covstrongdep} reports, for the case of strong serial dependence,  the empirical coverage of $95$\% Symmetric Confidence Intervals ($95$\%-SCI), for $\gamma_0$, $\theta_0$, $R_{0}^{(n)} (0.9)$ (RL)—the return level corresponding BM of size $n$ with a return period of $10n$ days, $Q_{0} (1 - 1 / n)$ (EQ)—the VaR with single daily probability $1 - 1 / n$, defined in formulas \eqref{eq:CI_MLE_gamma}, \eqref{eq:CI_MLE_RL} and \eqref{eq:CI_MLE_EXTQ}, respectively, which we will refer to as frequentist Symmetric (FS).  It also reports the empirical coverage of $95$\%-Symmetric-Bayesian Credible Intervals ($95$\%-SBCI), defined in \eqref{eq:BCI_gamma}, \eqref{eq:SBCI_EXT_IND}, \eqref{eq:SBCI_RL} and \eqref{eq:SBCI_EXT_VAR}, which we will refer to as Bayesian Symmetric (BS), and Asymmetric Bayesian Credible Intervals ($95$\%-SBCI and $95$\%-ABCI), defined in formulas \eqref{eq:credible_int_GEV}, \eqref{eq:ABCI_EXT_IND}, \eqref{eq:ABCI_RL},  and \eqref{eq:ABCI_EXT_VAR}, which we will refer to as Bayeisan Asymmetric (BA).  Finally, it also reports the empirical coverage $95$\% Asymmetric Confidence Intervals ($95$\%-ACI) for RL and EQ, defined in formulas \eqref{eq:ACI_RL} and \eqref{eq:ACI_EXTQ},  which we will refer to as frequentist Symmetric (FA). For simplicity, neglect asymptotic bias ($b_0=0$). Due to space constraints, the results obtained with the ARCH model are presented in the supplement; here, we provide a brief summary.

\begin{table}[t!]
	{\scriptsize
		\centering
		\caption{Mean squared error ratio (strong dependence)}
		\resizebox{\textwidth}{!}{\begin{tabular}{c c c|c c c c | c c c c | c c c c}
				\hline
				\multicolumn{3}{c}{} & \multicolumn{4}{c}{} & \multicolumn{8}{c}{Markov copula model ($\eta = 1.06$)} \\
				\multicolumn{3}{c}{} & \multicolumn{4}{c}{ARMAX ($\eta = 0.5$)} & \multicolumn{4}{c}{Exponential distr.} &  \multicolumn{4}{c}{Powerlaw distr.} \\
				\hline
				$n$ & $k$ & $m$ 
				& $\gamma_0$ & $\theta_0$ &  RL & EQ 
				& $\gamma_0$ & $\theta_0$ &  RL & EQ
				& $\gamma_0$ & $\theta_0$ &  RL & EQ \\
				\hline
				360 & 12 & 30 & 0.42 & 0.94 & 0.00 & 0.00 & 0.55 & 1.00 & 0.00 & 0.20 & 0.96 & 1.00 & 0.00 & 0.22 \\
				& 51 & 7 & 0.90 & 1.11 & 0.37 & 0.51 & 0.98 & 1.06 & 0.96 & 0.96 & 1.17 & 1.06 & 1.29 & 1.09 \\
				\hline
				1800 & 20 & 90 & 0.68 & 1.14 & 0.00 & 0.01 & 0.73 & 1.15 & 0.54 & 0.81 & 0.97 & 1.15 & 1.12 & 1.06 \\
				& 60 & 30 & 0.90 & 1.09 & 0.51 & 0.59 & 0.97 & 1.09 & 1.04 & 1.02 & 1.17 & 1.09 & 1.27 & 1.07 \\
				\hline
				3600 & 20 & 180 & 0.68 & 1.14 & 0.00 & 0.00 & 0.71 & 1.16 & 0.02 & 0.16 & 0.87 & 1.16 & 0.82 & 0.91 \\
				& 40 & 90 & 0.85 & 1.12 & 0.30 & 0.41 & 0.89 & 1.12 & 1.05 & 1.03 & 0.97 & 1.12 & 1.22 & 1.04 \\
				& 120 & 30 & 0.94 & 1.06 & 0.70 & 0.75 & 1.02 & 1.05 & 1.08 & 1.04 & 1.19 & 1.05 & 1.23 & 1.09 \\
				\hline
				10800 & 30 & 360 & 0.80 & 1.14 & 0.20 & 0.31 & 1.02 & 1.05 & 0.91 & 1.04 & 0.67 & 1.14 & 0.56 & 0.43 \\
				& 60 & 180 & 0.90 & 1.08 & 0.50 & 0.59 & 0.82 & 1.14 & 1.08 & 0.98 & 0.77 & 1.10 & 0.63 & 0.54 \\
				& 120 & 90 & 0.95 & 1.05 & 0.69 & 0.75 & 0.93 & 1.10 & 1.08 & 1.05 & 1.11 & 1.07 & 1.14 & 1.03 \\
				& 360 & 30 & 0.98 & 1.04 & 0.86 & 0.89 & 0.93 & 1.10 & 1.07 & 1.05 & 1.14 & 1.02 & 1.20 & 1.13 \\
				\hline
		\end{tabular}}
		\label{table:msestrongdep}
	}
	\footnotesize{
		Ratio of empirical mean squared errors of the posterior median versus the MLE for $\gamma_0$,  $\theta_0$,  $R_{0}^{(n)}(0.9)$ (RL), and $Q_{0}(1 - 1/n)$ (EQ).}
\end{table}

Encouragingly, the coverage of the 95\%-SBCI and 95\%-ABCI intervals for  $\gamma_0$ and the RL is already close to the nominal level at the smallest setting $n = 360$, $k = 12$, $m = 30$, across all time series models except for the ARCH model where good coverage is only achieved for $n \geq 1800$. As expected, performance improves with increasing of $n$. For EQ, coverage is initially lower than the nominal level with the smallest setting, but approaches the nominal level starting from the moderate setting $n = 1800$ with $k = 20$ and $m = 90$, again across all models.
For  $\theta_0$, performance varies by model. In the ARMAX setting, coverage begins to stabilize around the nominal level at moderate setting, while for the Markov copula models and the ARCH model acceptable coverage is only achieved at the largest setting $n = 10800$, with $k = 30$ and $m = 360$. 
The 95\%-SCI for $\gamma_0$ also performs well, matching nominal levels even at $n = 360$, and respectively $n = 1800$ for the ARCH model, except when the block size $m$ is much smaller than the number of blocks $k$. For $\theta_0$, 95\%-SCI performance again varies by model, with better results in the ARMAX case than in the Markov copula models and ARCH model, particularly when $k \gg m$. Overall, the 95\%-SCI for $\theta_0$ tend to outperform their Bayesian counterparts slightly, achieving around 91\% coverage at moderate setting $n = 1800$ with $k = 20$ and $m = 90$.
In contrast, the 95\%-SCI for RL and EQ consistently underperform relative to the 95\%-SBCI and 95\%-ABCI, especially in the ARMAX model, where nominal level is only reached at the largest setting $n = 10800$ with $k = 360$ and $m = 30$ and in the ARCH model, where coverage never exceeds the level $0.90$. The 95\%-ACI for RL and EQ generally outperform their symmetric counterparts and show comparable performance to Bayesian intervals, particularly in the ARMAX and ARCH case. However, for the Markov copula models, the Bayesian intervals retain a slight advantage in covering EQ. 

For brevity, results under weak serial dependence are reported in Table 4 and 5 of the Supplement, while here we only report a summary. 
Overall, all interval types—95\%-SBCI, 95\%-ABCI, 95\%-SCI, and 95\%-ACI—perform better under weak dependence. An exception is the ARCH model with weak dependence, where the coverage for $\gamma_0$ is significantly improved. Notably, both 95\%-ABCI and 95\%-SCI intervals for $\theta_0$ achieve near-nominal level even at the smallest setting for both ARMAX and Markov copula models.
As with strong dependence, 95\%-SCI  for RL and EQ perform significantly worse than 95\%-ACI, 95\%-SBCI, and 95\%-ABCI. Bayesian methods retain a slight edge over asymmetric intervals, particularly for EQ in the Markov copula and ARCH models.

Table \ref{table:msestrongdep} (Table 7 for the ARCH model) reports the ratio of empirical mean squared errors between the posterior median and the MLE under strong serial dependence. A ratio below one indicates that the Bayesian point estimator outperforms the MLE, and vice versa.
For the ARMAX and ARCH models and the Markov copula model with exponentially distributed marginals, the posterior median for $\gamma_0$ consistently outperforms the MLE, particularly with the smallest and moderate settings $n = 360,1800, 3600$. In contrast, with the Markov copula model and power-law marginals, the results are more nuanced: the posterior median performs better when $m > k$, whereas the MLE tends to perform better when $k > m$.
For $\theta_0$, the MLE generally has a slight advantage across all models. However, the most pronounced differences are seen for the RL and EQ estimates in the ARMAX and ARCH models, where the posterior median significantly outperforms the MLE, especially when $m \gg k$ and at smallest $m$. A similar, though less pronounced, pattern holds for the copula models with lighter tails.
With larger sample sizes and more balanced values of $k$ and $m$, the performance gap narrows, and in short-tailed cases, the MLE occasionally outperforms the posterior median, particularly for the RL.
The corresponding results under weak dependence, provided in Table 6 and 7 of the Supplement, largely mirror the patterns observed under strong dependence. One notable exception is in the estimation of $\theta_0$ for the ARMAX model when $\theta_0 = 1$, where the posterior median clearly outperforms the MLE. This improvement is attributed to the use of a prior $\Pi$ with a point mass of $p = 0.1$ at $\theta = 1$, which enhances estimation accuracy in this boundary case.

\begin{figure}[t!]
	\centering
	\includegraphics[page=1,width=.24\textwidth]{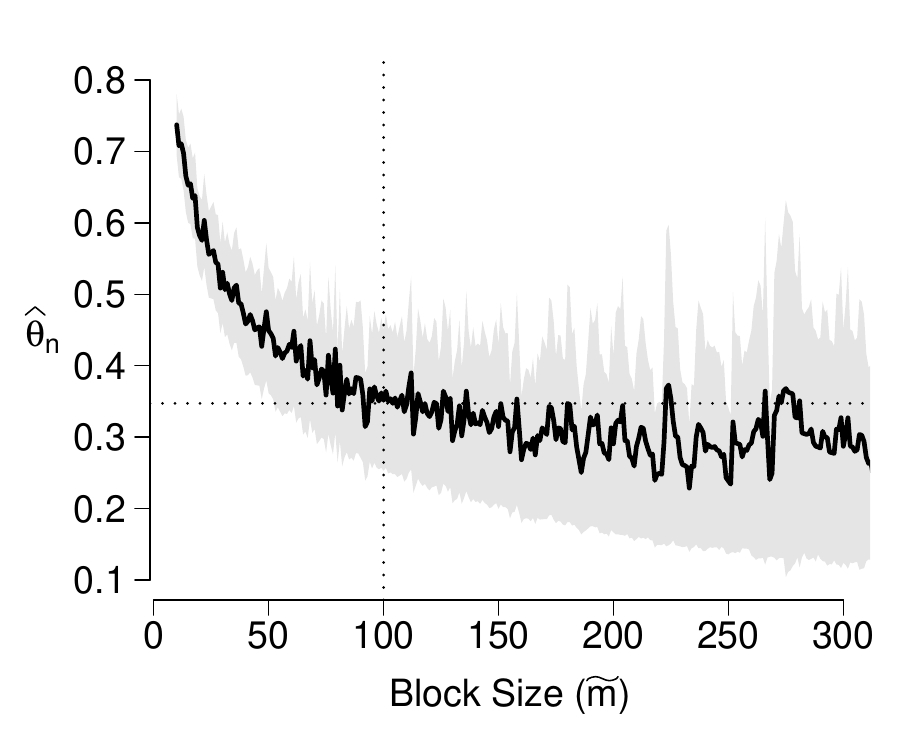}
	\includegraphics[page=2,width=.24\textwidth]{real_data_example.pdf}
	\includegraphics[page=3,width=.24\textwidth]{real_data_example.pdf}
	\includegraphics[page=4,width=.24\textwidth]{real_data_example.pdf}\\
	\includegraphics[page=5,width=.24\textwidth]{real_data_example.pdf}
	\includegraphics[page=6,width=.24\textwidth]{real_data_example.pdf}
	\includegraphics[page=7,width=.24\textwidth]{real_data_example.pdf}
	\includegraphics[page=8,width=.24\textwidth]{real_data_example.pdf}
	\caption{From top-left to top-right: extremal index estimate  $\widehat{\theta}$ (left) and extreme value index estimate $\widehat{\gamma}$ vs block size with 95\%-SCI in grey, 
		estimated autocorrelation function for $\widehat{m}_{i, m}$ (black lines) and $\widehat{m}_{i, m}^2$ (grey lines), with 95\%-SCI (dashed horizontal lines), and
		Q-Q plot of $F_{n}^{\widehat{\theta} m} (\widehat{m}_{i, m})$, for $i=1,\ldots,k$. From bottom-left to bottom-right: empirical posterior and prior density for $\gamma_0$, $\theta_{0}^{-1}$, and $Q_{0} (0.999)$ with corresponding medians (black squares) and 95\%-ABCI (black diamonds),
		and the daily return of the MSCI World Index with posterior median and $95$\%-ABCI of $Q_{0} (0.999)$ superimposed.}
	\label{fig:real_data_example1}
\end{figure}

%
%
\section{Real data example}
%
%
We analyze daily returns (negative log-returns) of the MSCI World Index over the period 2000–2024, yielding $n=5217$ observations. Our primary goal is to perform Bayesian inference on three key quantities: the tail index $\gamma_0$, the mean cluster size of extreme losses  $\theta_{0}^{-1}$, and the extreme VaR at level 0.999, denoted $Q_{0} (0.999)$.
Previous studies on such daily financial time series suggest that their marginal distributions are heavy-tailed, although it is typically assumed to possess a finite second moment, e.g. \citet[]{cont2001}. However, the condition $0<\gamma_0 < 0.5$, which would imply a finite second moment, has been questioned in the literature; see, for instance, \citet[][Ch. 7.6]{embrechts1997modelling}.

We adopt the Bayesian approach detailed in Section \ref{sec:Bayes_inf}, using data-driven priors for the location $\mu$ and scale $\sigma$ parameters as described in Section \ref{sec:simulations}. For $\gamma_0$, we employ the prior density $\pi_{\text{sh}}^{(n)}$, defined as a truncated Gaussian centered at 0.25 with variance 0.16, constrained to the interval $(-1/2, \sqrt{k})$. This prior allocates approximately 50\% of its mass to the interval (0,0.5), representing a heavy-tailed density with a finite second moment, while the remaining mass allows for either a light-tailed density or a heavy-tailed density with infinite variance.
The prior for $\theta_0$ is defined as a mixture of a uniform distribution on $(0,1)$ and the point mass  $0.1$ at 1, as introduced in Section \ref{sec:simulations}. Credible intervals for $\theta_{0}^{-1}$ are obtained via the transformation $x \rightarrow x^{-1}$ for $x> 0$, applied to the posterior of $\theta_0$. Since such a transformation is continuous and strictly monotone, then the coverage of the resulting intervals is also asymptotically correct.

Since the MLE and Bayesian asymptotics share the same conditions, we use $\widehat{\theta}_n$ and  $\widehat{\bfvartheta}_n$ to guide the choice of the big-blocks $\lm$ and $m$ and the small-block $l$. Following \cite{berghaus2018}, we compute the estimate  $\widehat{\theta}$ for $\lm = 10, \dots, 348$ and select the smallest $\lm$ for which the estimates appear stable. Based on the top-left panel of Figure \ref{fig:real_data_example1}, we choose $\lm = 104$, corresponding roughly to half-yearly maxima.
Rather than verifying Conditions \ref{cond:ts_cond}\ref{cond_alpha_mix1}–\ref{cond_beta_mix} in full detail, we adopt a pragmatic approach. In essence, Conditions \ref{cond:ts_cond}\ref{cond_alpha_mix1}–\ref{cond_alpha_mix2} ensure asymptotic independence across blocks, while Condition \ref{cond:ts_cond}\ref{cond_beta_mix} implies that block maxima behave like iid maxima over blocks of size $\lfloor\theta_0 \lm \rfloor$. We begin with $m = 104$ and $l = 0$, and compute autocorrelations of the empirical block maxima $\widehat{m}_{i,m}$ and $\widehat{m}_{i,m}^{2}$ for $i=1,\ldots,k$. We also generate Q–Q plots of the empirical distribution $F_{n}^{\widehat{\theta} m} (\widehat{m}_{i,m})$, $i = 1, \dots, k$ against the uniform one. These diagnostics are repeated for increasing values of $m$ and $l$ until autocorrelations lie within 95\%-SCI and the Q–Q plot shows no systematic deviation from the diagonal. Based on the top-third and top-right panels of Figure \ref{fig:real_data_example1}, we confirm the initial choice $m=104$ and  $l=0$, yielding $k = \lk = 50$ blocks.
To validate Condition \ref{cond:second_order}, we evaluate the estimate $\widehat{\gamma}$ over  $m = 104, \dots, 348$. Stability of $\widehat{\gamma}$ for $m\geq104$, as seen in the top-second panel of Figure \ref{fig:real_data_example1}, suggests that the condition holds as bias due to block size is negligible.

The bottom panels of Figure \ref{fig:real_data_example1} displays the prior (dotted lines) and the empirical posterior (solid lines) distributions for $\gamma_0$, $\theta_{0}^{-1}$, and $Q_{0} (\tau_E)$, based on 100,000 posterior draws. See Section E of the Supplement for details about the sampling scheme. Posterior median and 95\%-ABCI  are reported with black squares and diamonds, respectively.
Consistent with prior findings, the negative log-returns of the MSCI World Index are confirmed to follow a heavy-tailed distribution with finite first moment, as evidenced by the 95\%-ABCI for $\gamma_0$, which is $[0.08; 0.59]$ and excludes the value $1$. The inclusion of values between $[0.5; 0.59]$ suggests that the assumption of a finite second moment is not strongly supported by the data.
The posterior distribution of the mean cluster size $\theta_{0}^{-1}$ suggests substantial extremal dependence, as the 95\%-ABCI of $[2.2; 3.8]$ implies that exceedances over high thresholds are expected to occur in clusters of two to approximately four. The posterior for the extreme VaR $Q_{0} (0.999)$ is strongly right-skewed, with median 0.065 and $[0.048; 0.111]$ as 95\%-ABCI. The bottom-right panel of Figure \ref{fig:real_data_example1} displays the daily returns of the MSCI World Index, with the posterior median and 95\%-ABCI superimposed. These results are informative, as the interval contains the largest observed return and captures 0.17\% of the observed returns, as expected, offering a plausible range for the intensity of future extreme events.

\section*{Acknowledgement}
Simone Padoan is supported by the Bocconi Institute for Data Science and
Analytics (BIDSA) and MUR–PRIN
Bando 2022–prot. 20229PFAX5, financed by the European Union - Next Generation
EU, Mission 4 Component 2 CUP J53D23004260001.

\section*{Supplementary material}
The Supplementary Material document contains a detailed discussion of important frequentist properties of the empirical loglikelihood process relative to the GEV model, revisited asymptotic results obtained under weaker conditions, all the proofs, further numerical details and results concerning the simulation study.

\bibliographystyle{chicago} 
\bibliography{references}

\end{document}